\documentclass[a4paper,12pt]{amsart}
\usepackage{amsmath, amsfonts,amsthm,times}
\usepackage{mathrsfs}
\usepackage[active]{srcltx}
\usepackage{amssymb}
\usepackage{ifthen}
\usepackage[usenames]{color}
\usepackage{graphicx}
\nonstopmode \numberwithin{equation}{section}
\newtheorem{thm}{Theorem}[section]

\newtheorem{cor}{Corollary}[section]
\newtheorem{lem}{Lemma}[section]
\newtheorem{prop}{Proposition}[section]

\newtheorem{claim}{Claim}[section]
\newtheorem{subclaim}{Subclaim}
\newtheorem{prob}{Problem}[section]
\newtheorem{case}{Case}[section]
\newtheorem*{mysolution}{Solution}
\newtheorem{step}{Step}[section]
\theoremstyle{definition}
\newtheorem{eg}{Example}[section]
\newtheorem{defn}{Definition}[section]
\newtheorem{examp}{Example}[section]

\newtheorem{ques}{Question}[section]
\newtheorem{rem}{Remark}[section]
\newcounter {own}
\def\theown {\thesection       .\arabic{own}}

\newenvironment{pf}[1][]{%
 \vskip 3mm
 \noindent
 \ifthenelse{\equal{#1}{}}%
  {{\slshape Proof. }}%
  {{\slshape #1.} }%
 }%
{\qed\bigskip}

\newcounter{alphabet}

\newenvironment{Thm}[1][]{\refstepcounter{alphabet}%
\bigskip%
\noindent%
{\bf Theorem \Alph{alphabet}}%
\ifthenelse{\equal{#1}{}}{}{ (#1)}%
{\bf .} \itshape}{\vskip 8pt}

\makeatletter

\newenvironment{Lem}[1][]{\refstepcounter{alphabet}%
\bigskip%
\noindent%
{\bf Lemma \Alph{alphabet}}%
{\bf .} \itshape}{\vskip 8pt}

\def\be{\begin{equation}}
\def\ee{\end{equation}}

\newcommand{\ben}{\begin{enumerate}}
\newcommand{\een}{\end{enumerate}}

\newcommand{\blem}{\begin{lem}}
\newcommand{\elem}{\end{lem}}
\newcommand{\bthm}{\begin{thm}}
\newcommand{\ethm}{\end{thm}}
\newcommand{\bcor}{\begin{cor}}
\newcommand{\ecor}{\end{cor}}
\newcommand{\beg}{\begin{examp}}
\newcommand{\eeg}{\end{examp}}
\newcommand{\begs}{\begin{examples}}
\newcommand{\eegs}{\end{examples}}
\newcommand{\bdefe}{\begin{defn}}
\newcommand{\edefe}{\end{defn}}
\newcommand{\bprob}{\begin{prob}}
\newcommand{\eprob}{\end{prob}}
\newcommand{\bques}{\begin{ques}}
\newcommand{\eques}{\end{ques}}
\newcommand{\bei}{\begin{itemize}}
\newcommand{\eei}{\end{itemize}}
\newcommand{\bcl}{\begin{claim}}
\newcommand{\ecl}{\end{claim}}
\newcommand{\bscl}{\begin{subclaim}}
\newcommand{\escl}{\end{subclaim}}
\newcommand{\bca}{\begin{case}}
\newcommand{\eca}{\end{case}}
\newcommand{\bstep}{\begin{step}}
\newcommand{\estep}{\end{step}}
\newcommand{\bsol}{\begin{mysolution}}
\newcommand{\esol}{\end{mysolution}}
\newcommand{\bcon}{\begin{conj}}
\newcommand{\econ}{\end{conj}}
\newcommand{\bcons}{\begin{conjs}}
\newcommand{\econs}{\end{conjs}}
\newcommand{\bprop}{\begin{prop}}
\newcommand{\eprop}{\end{prop}}
\newcommand{\br}{\begin{rem}}
\newcommand{\er}{\end{rem}}
\newcommand{\brs}{\begin{rems}}
\newcommand{\ers}{\end{rems}}
\newcommand{\bo}{\begin{obser}}
\newcommand{\eo}{\end{obser}}
\newcommand{\bos}{\begin{obsers}}
\newcommand{\eos}{\end{obsers}}
\newcommand{\bpf}{\begin{pf}}
\newcommand{\epf}{\end{pf}}
\newcommand{\ba}{\begin{array}}
\newcommand{\ea}{\end{array}}
\newcommand{\beq}{\begin{eqnarray}}
\newcommand{\beqq}{\begin{eqnarray*}}
\newcommand{\eeq}{\end{eqnarray}}
\newcommand{\eeqq}{\end{eqnarray*}}

\begin{document}
\bibliographystyle{amsplain}

\title{Isoperimetric type  inequalities for mappings induced by   weighted Laplace differential operators}

\author{Jiaolong Chen}
\address{Jiaolong Chen,
MOE-LCSM, School of Mathematics and Statistics, Hunan Normal University, Changsha, Hunan 410081, P. R. China}
\email{jiaolongchen@sina.com}

\author{Shaolin Chen$^{\mathscr{*}}$}
\address{Shaolin Chen,
College of Mathematics and Statistics, Hengyang Normal University, Hengyang, Hunan 421002, P. R. China}
\email{mathechen@126.com}

\author{Manzi Huang}
\address{Manzi Huang,
MOE-LCSM, School of Mathematics and Statistics, Hunan Normal University, Changsha, Hunan 410081, P. R. China} \email{mzhuang@hunnu.edu.cn}

\author{HUAQING ZHENG}
\address{Huaqing Zheng,
MOE-LCSM, School of Mathematics and Statistics, Hunan Normal University, Changsha, Hunan 410081, P. R. China}
 \email{huaqingzheng2022@163.com}


\keywords{Isoperimetric type  inequality, Poisson type integral,   Hardy type  space,   Bergman type  space, $(K,K')$-elliptic mapping.\\
${}^{\mathbf{*}}$ Corresponding author}

\begin{abstract}
The main purpose of this paper is to establish some  isoperimetric type  inequalities for  mappings induced by  the weighted Laplace differential operators.
The obtained results of this paper provide  improvements and
extensions of the corresponding known results.
\end{abstract}

\maketitle \pagestyle{myheadings} \markboth{Jiaolong Chen, Shaolin Chen, Manzi Huang and Huaqing Zheng}{Isoperimetric type  inequalities for mappings induced by  Poisson type integrals}

\section{Introduction and main results }\label{sec-1}
Denote by  $\mathbb{C}$  the complex plane. Let $\mathbb{D}=\{z\in \mathbb{C}:|z|<1\}$ be the unit disk, and let $\mathbb{T}:=\partial\mathbb{D}$ be the unit circle.
For a domain $\Omega\subset\mathbb{C}$,
we use $\mathcal{C}^n(\Omega)$ to denote the set of
all $n$-times continuously differentiable complex-valued functions  in $\Omega$,
where $n\in \mathbb{N}_{0}=\{0\}\cup \mathbb{N}$, and $\mathbb{N}$ denotes the set of all positive integers.
In particular, let $\mathcal{C}(\Omega):=\mathcal{C}^{0}(\Omega)$. For $z=x+iy\in\mathbb{C}$, the complex formal differential
operators are defined by $\partial/\partial z=1/2(\partial/\partial x-i\partial/\partial y)$ and
$\partial/\partial \overline{z}=1/2(\partial/\partial x+i\partial/\partial y)$, where $x,~y\in\mathbb{R}$.
For
$\theta\in[0,2\pi]$, the  directional derivative of  a
mapping ( or complex-valued function) $f\in\mathcal{C}^n(\mathbb{D})$ at $z\in\mathbb{D}$ is defined
by

$$\partial_{\theta}f(z):=\lim_{\rho\rightarrow0^{+}}\frac{f(z+\rho e^{i\theta})-f(z)}{\rho}=\partial_{z}f(z)e^{i\theta}
+\partial_{\overline{z}}f(z)e^{-i\theta},$$ where $\partial_{z}f:=\partial
f/\partial z,$ $\partial_{\overline{z}}f:=\partial f/\partial \overline{z}$
and $\rho$ is a positive real number such that $z+\rho
e^{i\theta}\in\mathbb{D}$. Let

$$\|D_{f}(z)\|:=\sup\{|\partial_{\theta}f(z)|:\; \theta\in[0,2\pi]\}
~\mbox{and}
~l\big(D_{ f}(z)\big):=\inf\{|\partial_{\theta}f(z)|:\; \theta\in[0,2\pi]\}.$$
Then

\begin{equation}\label{eq-1.1}
\|D_{f}(z)\|=|\partial_{z}f(z)|+|\partial_{\overline{z}}f(z)|\;\;\mbox{and}\;\; l\big(D_{ f}(z)\big)=\big ||\partial_{z}f(z)|-|\partial_{\overline{z}}f(z)|\big |,
\end{equation}



\subsection{The weighted Laplace differential operators}
For $z\in\mathbb{D}$, let $$\sigma_{\alpha}(z)=(1-|z|^{2})^{\alpha}$$
be the standard weight in $\mathbb{D}$, where $\alpha\in\mathbb{R}$. Let us recall the
differential operator $\Delta_{\alpha}$ related to the  standard weight $\sigma_{\alpha}$ as follows:
\be\label{C-y}\Delta_{\alpha}:=\frac{1}{4}{\rm div}(\sigma_{\alpha}^{-1}\nabla)-\frac{\alpha^{2}}{4}\sigma_{\alpha+1}^{-1},\ee
where the symbols $\nabla$ and ${\rm div}$ denote the gradient and divergence, respectively.
The formula (\ref{C-y}) points at a relation to the so-called conductivity equations considered by
Astala and P\"aiv\"arinta \cite{A-P} (see also \cite{ao}). Of particular interest to our analysis is the {\it Dirichlet problem}:


\begin{align}\label{eq-1.2}
\begin{cases}
\displaystyle \Delta_{\alpha}f =0 \;\;&\text{in} \;\; \mathbb{D}, \\
\displaystyle f=F \;\;&\text{on}\;  \;\mathbb{T}.
\end{cases}
\end{align}
Here, the boundary data $F\in \mathcal{D}'(\mathbb{T})$ is a \emph{distribution}
on $\mathbb{T}$, and the boundary condition in (\ref{eq-1.2}) is interpreted
in the distributional sense that $f_{r} \rightarrow F $
in $\mathcal{D}'(\mathbb{T})$ as $r\rightarrow 1^{-}$ (see Section 2 below for the details), where $f_{r}( e^{i\theta}):=f(re^{i\theta})$ for $\theta\in[0,2\pi]$.
It is shown  that, for $\alpha>-1$,
if a mapping $f\in \mathcal{C}^2(\mathbb{D})$ satisfies $\Delta_{\alpha}f=0$ in $\mathbb{D}$
and
 $f_{r} \rightarrow F $
in $\mathcal{D}'(\mathbb{T})$ as $r\rightarrow 1^{-}$, then $f$ can
be expressed by the following Poisson type integral:
\be\label{eq-1.3}
f =\mathcal{K}_{\alpha}[F]
\ee
in $ \mathbb{D}$,
where
\be\label{eq-1.4}
\mathcal{K}_{\alpha}[F](z)=\frac{1}{2\pi}\int_{0}^{2\pi}\mathcal{K}_{\alpha}(ze^{-it})F(e^{it})\,dt,
\quad \mathcal{K}_{\alpha}(z)= \mathcal{C}_{\alpha}\frac{(1-| z|^2)^{\alpha+1}}{| 1-z|^{\alpha+2}},
\ee
 \begin{equation}\label{eq-1.5}
 \mathcal{C}_{\alpha}=\frac{  \Gamma ^2(1+\frac{\alpha}{2} )}{\Gamma(1+\alpha)},
 \end{equation} and
 $\Gamma$  is the  Gamma function (see \cite{bh,ao}).
If $\alpha\leq-1$,  then $f\equiv0$ in $\mathbb{D}$ (see \cite[Theorem 2.3]{ao}).
 In particular, if
$\alpha=0$ in (\ref{eq-1.3}), then $f$  is harmonic in  $\mathbb{D}$.
It is well known that every harmonic mapping $f$ defined in  $\mathbb{D}$ admits the canonical
decomposition $f=h+\overline{g}$, where $h$ and $g$ are analytic in $\mathbb{D}$.
This representation is unique up to an additive constant (cf. \cite[p.7]{pd}).
We refer to \cite{bh, kmm,L-W, A-2020,ao} for basic properties of
the mappings defined in (\ref{eq-1.3}). For higher dimensional
cases, see \cite{LP,LP-2009,liu}.

\subsection{Hardy type spaces and  Bergman type spaces }
For $p\in (0, \infty]$, the {\it  Hardy type   space}
$\mathcal{H}_{\mathcal{G}}^p(\mathbb{D})$ consists of all
measure complex-valued functions $f$ defined in $\mathbb{D} $ such that
$H_p(f)<\infty$,
where $H_p(f)=\sup_{r\in[0,1)}M_p(r,f)$, and
$$ M_p(r,f)=\begin{cases}
\displaystyle \left(\frac{1}{2\pi}\int_{0}^{2\pi}| f(re^{i\theta})| ^p\,d\theta\right)^{\frac{1}{p}},\;\;&\text{if} \;\;p\in(0,\infty), \\
\displaystyle \max_{ \theta\in [0,2\pi]} \big| f(re^{i\theta})\big| ,\;\;&\text{if}\;\; p=\infty
\end{cases}
$$
(cf. \cite{cpr}).
 The norm in $\mathcal{H}_{\mathcal{G}}^p(\mathbb{D})$ is denoted by $\|\cdot\| _{\mathcal{H}_{\mathcal{G}}^p(\mathbb{D})}$ which is defined by
$$\|f\| _{\mathcal{H}_{\mathcal{G}}^p(\mathbb{D})}=H_p(f)$$ for every $f\in \mathcal{H}_{\mathcal{G}}^p(\mathbb{D})$.

For $p\in (0, \infty]$, the {\it Bergman type   space}
$\mathcal{B}_{\mathcal{G}}^p(\mathbb{D})$ consists of all
measure complex-valued functions $f$ defined in $\mathbb{D} $ such that
$B_p(f) <\infty$,
where
$$B_p(f)=\begin{cases}
\displaystyle \left(\int_{\mathbb{D}}| f(z)|^p\,d\sigma(z)\right)^{\frac{1}{p}},\;\;&\text{if} \;\;p\in(0,\infty), \\
\displaystyle \text{esssup}_{z\in \mathbb{D}}|f(z)|,\;\;&\text{if}\;\; p=\infty,
\end{cases}
$$
and $d\sigma=dxdy/\pi$ is the normalized Lebesgue measure in $\mathbb{D}$ (cf. \cite[p.1]{hhb}).
  The norm in $\mathcal{B}_{\mathcal{G}}^p(\mathbb{D})$ is denoted by $\|\cdot\| _{\mathcal{B}_{\mathcal{G}}^p(\mathbb{D})}$ which is defined by
 $$
 \|f\| _{\mathcal{B}_{\mathcal{G}}^p(\mathbb{D})}=B_p(f)
 $$ for every $f\in \mathcal{B}_{\mathcal{G}}^p(\mathbb{D})$.

The classical Hardy space $\mathcal{H}^p(\mathbb{D})$ and the classical Bergman space $\mathcal{B}^p(\mathbb{D})$,
i.e., all their elements are analytic,
are the subspaces of $\mathcal{H}_{\mathcal{G}}^p(\mathbb{D})$ and $\mathcal{B}_{\mathcal{G}}^p(\mathbb{D})$,
respectively.
If $f$ is analytic in $\mathbb{D}$,
then
$$
\|f\| _{\mathcal{H}_{\mathcal{G}}^p(\mathbb{D})}=\|f\| _{\mathcal{H}^p(\mathbb{D})}
\;\;\text{and}\;\;
\|f\| _{\mathcal{B}_{\mathcal{G}}^p(\mathbb{D})}=\|f\| _{\mathcal{B}^p(\mathbb{D})}
$$
(cf. \cite[Section 1.1]{dp} and \cite[Section 1.1]{hhb}).
Obviously, for $0<p_{1}<p_{2}\leq \infty$, we have the following inclusions:
\be\label{eq-1.6}
\mathcal{H}_{\mathcal{G}}^{p_{2}}(\mathbb{D})
 \subset\mathcal{B}_{\mathcal{G}}^{p_{2}}(\mathbb{D}),
\ee
$$\mathcal{H}_{\mathcal{G}}^{\infty}(\mathbb{D})
 \subset\mathcal{H}_{\mathcal{G}}^{p_{2}}(\mathbb{D})
 \subset
\mathcal{H}_{\mathcal{G}}^{p_{1}}(\mathbb{D})
\;\;\text{and}\; \;
\mathcal{B}_{\mathcal{G}}^{\infty}(\mathbb{D})
\subset\mathcal{B}_{\mathcal{G}}^{p_{2}}(\mathbb{D})
\subset
\mathcal{B}_{\mathcal{G}}^{p_{1}}(\mathbb{D}).
$$

For $p\in (0, \infty)$, we use $\mathcal{L}^p(\mathbb{T})$ to denote the space of
measure complex-valued functions $F$ defined on $\mathbb{T} $ such that
$$\int_{0}^{2\pi}| F(e^{i\theta})| ^p\,d\theta<\infty $$
(cf. \cite[p.65]{rudin1987}).  The norm in $\mathcal{L}^p(\mathbb{T})$ is denoted by $ \|\cdot\| _{\mathcal{L}^p(\mathbb{T})}$ which is defined by
\begin{equation}\label{eq-1.7}
 \|F\| _{\mathcal{L}^p(\mathbb{T})}= \left(\frac{1}{2\pi}\int_{0}^{2\pi}| F(e^{i\theta})| ^p\,d\theta\right)^{\frac{1}{p}}.
\end{equation}

When $p=\infty$, $\mathcal{L}^{\infty}(\mathbb{T})$ denotes the space of
measure   complex-valued functions $F$ defined on $\mathbb{T} $ with the norm
 $$
 \|F\| _{\mathcal{L}^{\infty}(\mathbb{T})}=\text{esssup} _{\theta\in [0,2\pi]}  \big| F(e^{i\theta})\big|
  $$
(cf. \cite[p.66]{rudin1987}).

The following known results will be used
later on.

\begin{Thm}\label{Thm-A}{\rm \bf (\cite[Theorem 5.1.8]{pav})}
If $f\in \mathcal{H}^{p}(\mathbb{D})$ for some $p\in(0,\infty]$,
then the radial limit $F(e^{i\theta})=\lim_{r\rightarrow1-} f(re^{i\theta})$ exists  a.e. $($i.e., almost everywhere$)$ on $\mathbb{T}$, and
 $$
 \|F \|_{ \mathcal{L}^{p}(\mathbb{T})}
 =\|f\|_{ \mathcal{H}^{p}(\mathbb{D})}
.
 $$
\end{Thm}

\begin{Thm}\label{Thm-B}{\rm \bf (\cite[Theorem 5.2.4]{pav})}
If $f\in \mathcal{H}^{1}(\mathbb{D})$ and the radial limit $F(e^{i\theta})=\lim_{r\rightarrow1^{-}} f(re^{i\theta})$  is
equal to a function of bounded variation a.e. on $\mathbb{T}$,
then $f$ has an absolutely continuous extension to $\overline{\mathbb{D}}=\{z\in \mathbb{C}:|z|\leq 1\}$.
\end{Thm}

 \begin{Thm}\label{Thm-C}{\rm \bf (\cite[Theorem 3.12]{dp})}
Let $f $ map $\mathbb{D}$ conformally onto a domain bounded by a Jordan curve $\gamma$.
Then $\gamma$ is rectifiable if and only if $f'\in \mathcal{H}^{1}(\mathbb{D})$.
\end{Thm}

Here, a {\it Jordan curve} (or a {\it simple closed curve}) $\gamma$ is the image of a continuous complex-valued function
$\gamma=\gamma(t)$ ($t\in[0,2\pi]$) such that
$\gamma(0)=\gamma(2\pi)$ and $\gamma(t_{1})\not=\gamma(t_{2})$ for $0\leq t_{1}\not=t_{2}<2\pi$.

It is well known that if $F$ is absolutely continuous on $\mathbb{T}$, then  the derivative
$$\dot F(e^{i\theta}):=\frac{d}{d\theta}F(e^{i\theta})$$
exists a.e. on $\mathbb{T}$ and $\dot F\in \mathcal{L}^{1}(\mathbb{T})$ (cf. \cite[Theorem 7.20]{rudin1987}).
Further, we have the following result.

 \begin{Thm}\label{Thm-D}{\rm \bf (\cite[Theorem 3.11]{dp})}
A function $f$ analytic in $\mathbb{D}$ is continuous in $\overline{\mathbb{D}}$      and
absolutely continuous on $\mathbb{T}$ if and only if  $f'\in \mathcal{H}^{1}(\mathbb{D})$.
If $f'\in \mathcal{H}^{1}(\mathbb{D})$, then
$$
\frac{d}{d\theta}F(e^{i\theta})= \lim_{r\rightarrow1^{^{-}}}\frac{\partial}{\partial\theta}f(re^{i\theta})
$$
a.e. on $\mathbb{T}$, where $F(e^{i\theta})=\lim_{r\rightarrow1-} f(re^{i\theta})$.
\end{Thm}

 \begin{Thm}\label{Thm-E}{\rm \bf (\cite[Theorem 4.1]{dp})}
If $\varphi:\mathbb{D}\rightarrow \mathbb{R}$ is a harmonic function and $\varphi\in \mathcal{H}_{\mathcal{G} }^{p}(\mathbb{D})$ for some $p\in(1,\infty)$,
then its harmonic conjugate $\varphi^*$ is also of class $\mathcal{H}_{\mathcal{G} }^{p}(\mathbb{D})$,
 where $\varphi^*(0)=0$.
Furthermore, there is a constant $A_{p}$, depending only on $p$, such that
$$
\|\varphi^*\|_{\mathcal{H}_{\mathcal{G}}^{p}(\mathbb{D})}\leq A_{p}\|\varphi\|_{\mathcal{H}_{\mathcal{G}}^{p}(\mathbb{D})}.
$$
\end{Thm}
\subsection{Elliptic Mappings}
A mapping $f:\Omega\rightarrow\mathbb{C}$ is said to be \emph{absolutely continuous on lines, $\mathcal{ACL}$} in brief, in the domain $\Omega$ if for every closed rectangle $R\subset\Omega$ with sides parallel to the axes $x$ and $y$, respectively, $f$ is absolutely continuous on almost every horizontal line and almost every vertical line in $R$. Such a mapping has, of course, partial derivatives $f_{x}$ and $f_{y}$ a.e. in $\Omega$. Moreover, we say that $f\in \mathcal{ACL}^2$ if its partial derivatives are locally Lebesgue square
 integrable in $\Omega$.

A sense-preserving and continuous mapping $f:\mathbb{D}\rightarrow\mathbb{C}$ is said to be a $(K,K')$-\emph{elliptic mapping} (or $(K,K')$-\emph{quasiregular mapping}) if\\
(1) $f$ is $\mathcal{ACL}^2$ in $\mathbb{D}$ and $J_{f}(z)\neq 0$ a.e. in $\mathbb{D}$, where $J_{f}$ denotes the Jacobian of $f$, which is given by
$$
J_{f}(z)=|\partial_{z}f(z) |^2-|\partial_{\overline{z}}f(z) |^2 =\|D_{f}(z) \| l \big(D_{f}(z)  \big);
$$
(2) there are constants $K\geq 1$ and $K'\geq 0$ such that
$$\| D_{f}(z)\|^2\leq KJ_{f}(z)+K'$$
 a.e. in $\mathbb{D}$.

In particular, if $K'\equiv 0$, then a $(K,K')$-elliptic mapping is said to be \emph{K-quasiregular}.
We refer to \cite{iwan, rick} for basic properties of $K$-quasiregular mappings.
It is well known that every quasiregular mapping is an elliptic mapping. However, the inverse of this statement is not true (see \cite{csw}). We refer to \cite{clsw, fs, kal} for more details of elliptic mappings.

\subsection{Isoperimetric type inequalities}
For a Jordan domain, we mean a domain whose boundary is a Jordan curve.
Let $\Omega\subset\mathbb{C}$ be a Jordan domain with  the rectifiable boundary $\partial\Omega$.
We use $ \mathcal{A}(\Omega)$ and  $ \ell(\partial\Omega)$ to denote the area of $\Omega$ and the length of  $\partial\Omega$, respectively.
Then there holds the isoperimetric inequality
 \be\label{eq-1.8}
\mathcal{A}(\Omega) \leq \frac{ \ell^{2}(\partial\Omega)}{4\pi},
\ee
 and equality occurs if and only if $\Omega$ is a disk (cf. \cite[p.25]{dp1983}).
Carleman \cite{carleman} gave a very beautiful proof of  the above isoperimetric inequality \eqref{eq-1.8} on minimal surfaces
 by reducing it to a theorem on analytic functions in $\mathbb{D}$.

Let $f$ be a conformal mapping from $\mathbb{D}$
onto a  Jordan domain $\Omega$ with  the   rectifiable boundary  $\partial\Omega$.
On the one hand,  we  obtain from Theorem C that
 $f'\in \mathcal{H}^{1}(\mathbb{D})$.
On the other hand, the Carath\'{e}odory extension theorem (cf. \cite[p.12]{dp1983})
asserts that $f$ can be extended to a homeomorphism on $\overline{\mathbb{D}}$.
Set
$$
F(e^{i\theta}) =\lim_{r\rightarrow1^{-}} f(re^{i\theta}).
$$
Then $F$ is a parametrization of $\partial\Omega$ (cf. \cite[p.44]{dp}),
and continuous  on $\mathbb{T}$.
It follows from \cite[Proposition 4.2]{zhu} that $f=\mathcal{K}_{0}[F]$.
Since $\partial\Omega$ is rectifiable, we see that $F$ is a mapping of bounded variation on $\mathbb{T}$.
This, together with Theorem B, implies  that $F$ is absolutely continuous on $\mathbb{T}$, and  $\dot F\in \mathcal{L}^{1}(\mathbb{T})$ (cf. \cite[Theorem 7.20]{rudin1987}),
where
$
\dot{F}(e^{i\theta})= \frac{d}{d\theta}F(e^{i\theta}).
$
Set $$\dot{f}(re^{i\theta}):= \frac{\partial}{\partial\theta}f(re^{i\theta}).$$
Then Theorem D implies that
\be\label{eq-1.9}
\dot{F}(e^{i\theta})= \lim_{r\rightarrow1^{^{-}}}\dot{f}(re^{i\theta})
\ee
  a.e.    on $\mathbb{T}$.
 Since Theorem C guarantees that $f'\in \mathcal{H}^{1}(\mathbb{D})$,
 we know from  Theorem  A  that
$$
 \| f'\|_{ \mathcal{H}^{1}(\mathbb{D}) }
 = \| \dot{f} \|_{ \mathcal{H}^{1}(\mathbb{D}) }
 =\Big\|\lim_{r\rightarrow1^{^{-}}}\dot{f}(re^{i\theta}) \Big\|_{\mathcal{L}^{1}(\mathbb{T}) }.
$$
This, together with  \eqref{eq-1.9}, yields that
$$
 \| f'\|_{ \mathcal{H}^{1}(\mathbb{D}) }
 =\Big\|\lim_{r\rightarrow1^{^{-}}}\dot{f}(re^{i\theta}) \Big\|_{\mathcal{L}^{1}(\mathbb{T}) }
 =\| \dot{F} \|_{\mathcal{L}^{1}(\mathbb{T}) }.
$$

In conclusion, we have the following proposition.

\begin{prop}\label{prop-1.1}
Suppose that $f$ is a conformal mapping of $\mathbb{D}$ onto a Jordan domain $\Omega$ with the rectifiable boundary, and $F$ denotes its boundary function on $\mathbb{T}$. Then the following statements are true.
\begin{enumerate}
\item[$(i)$]
$F $ is absolutely continuous on $\mathbb{T}$;
\item[$(ii)$]
$\dot{F} \in \mathcal{L}^{1}(\mathbb{T})$;
\item[$(iii)$]
$\| f'\|_{ \mathcal{H}^{1}(\mathbb{D}) }
=\|\dot{F}\|_{\mathcal{L}^{1}(\mathbb{T})}$.
\end{enumerate}
\end{prop}

By using the notations as above, \eqref{eq-1.8} can be reformulated as follows.

\begin{Thm}\label{Thm-F}{\rm \bf (\cite[p.82]{pav})}
Let $f$ be a conformal mapping from $\mathbb{D}$ onto a Jordan domain $\Omega$ with rectifiable boundary $\partial\Omega$.
Then
\be\label{eq-1.10}
\| f'\|_{\mathcal{B}^2(\mathbb{D})}
\leq \| f'\|_{ \mathcal{H}^{1}(\mathbb{D}) }
=\|\dot{F}\|_{\mathcal{L}^{1}(\mathbb{T})},
\ee
where $F(e^{i\theta}) =\lim_{r\rightarrow1^{-}} f(re^{i\theta})$.
\end{Thm}


By using a similar approach as in \cite{carleman}, Strebel generalized \eqref{eq-1.10}
 into the following  isoperimetric type inequalities:

\begin{Thm}\label{Thm-G}$($\cite[Theorem 19.9]{strebel}$)$ Let $f$ be an analytic function in $\mathbb{D}$  and
$f\in\mathcal{H}^{p}(\mathbb{D})$ with $p\in(0,\infty)$.
Then
\be\label{eq-1.11}
\| f\|_{\mathcal{B}^{2p}(\mathbb{D})}
\leq
\| f\|_{\mathcal{H}^{p}(\mathbb{D})}
= \|F\|_{\mathcal{L}^{p}(\mathbb{T})} ,
\ee
where  $F(e^{i\theta}) =\lim_{r\rightarrow1^{-}} f(re^{i\theta})$.
\end{Thm}
 Let $g'=f$ in Theorem G, and let $G(e^{i\theta})= \lim_{r\rightarrow1^{-}}g(re^{i\theta})$.
If $g$   is a conformal mapping from $\mathbb{D}$ onto a Jordan domain  with rectifiable boundary,
then  Theorems A and D ensure  that
$$
\|\dot{G}\|_{\mathcal{L}^{p}(\mathbb{T})}
=  \Big\|\lim_{r\rightarrow1^{^{-}}}\dot{g}(re^{i\theta}) \Big\|_{\mathcal{L}^{p}(\mathbb{T}) }
=\|F\|_{\mathcal{L}^{p}(\mathbb{T})}
=\| g'\|_{\mathcal{H}^{p}(\mathbb{D}) }.
$$
Then \eqref{eq-1.11} is equivalent to
\be\label{eq-1.12}
\| g'\|_{\mathcal{B}^{2p}(\mathbb{D})}
\leq \| g'\|_{ \mathcal{H}^{p}(\mathbb{D}) }
=\|\dot{G}\|_{\mathcal{L}^{p}(\mathbb{T})}.
\ee

Notice that $\| g'\|_{\mathcal{B}^{ p}(\mathbb{D})}\leq \| g'\|_{\mathcal{B}^{2p}(\mathbb{D})}$.
Therefore, \eqref{eq-1.12} implies the following isoperimetric type inequalities:
\be\label{eq-1.13}
\| g'\|_{\mathcal{B}^{ p}(\mathbb{D})}
\leq\| g'\|_{\mathcal{B}^{2p}(\mathbb{D})}
\leq \| g'\|_{ \mathcal{H}^{p}(\mathbb{D}) }
=\|\dot{G}\|_{\mathcal{L}^{p}(\mathbb{T})}.
\ee

Naturally, one will ask whether there are any results similar to those in \eqref{eq-1.12} or \eqref{eq-1.13} for mappings $f=\mathcal{K}_{\alpha}[F]$ which are defined in \eqref{eq-1.3}.
 In general, these mappings are not analytic. So we are going to replace the
 derivatives $f'$ in Theorem F by their partial derivatives:
$ \partial_{\theta}f$, $\partial_{r}f,$ $\partial_{z}f$ and $\partial_{\overline{z}}f$,
where  $\partial_{\theta}:=\partial /\partial \theta$ and
 $\partial_{r}:=\partial/\partial r$.
 Then we can state the problem as follows. In particular,  it is an open problem for $\alpha=0$ (see \cite[Problem 1]{zjf}).

 \bprob\label{prob-1.1}
 Suppose that
 \begin{enumerate}
   \item\label{prob-1.1-1}
   $F $ is an absolutely continuous function on $\mathbb{T}$
   and $\dot F\in \mathcal{L}^{p}(\mathbb{T})$ with $p \in [1,\infty]$;
   \item \label{prob-1.1-2}
   $f=\mathcal{K}_{\alpha}[F]$ in $\mathbb{D}$ with $\alpha\in(-1,\infty)$.
 \end{enumerate}
Do the partial derivatives $\partial_{\theta}f$, $\partial_{r}f$, $\partial_{z}f$ and $\partial_{\overline{z}}f$ of $f$ satisfy certain isoperimetric type inequalities similar to those in \eqref{eq-1.12} or \eqref{eq-1.13}?
\eprob

By Proposition \ref{prop-1.1}, we see that the assumptions in Problem \ref{prob-1.1}\eqref{prob-1.1-1} are natural.
 Hang et al. \cite[Theorem 1.1]{hang} obtained some isoperimetric type inequalities for harmonic functions in the upper half-space in $\mathbb{R}^{n}$.
Recently, Zhu \cite{zjf} and the second author et al. \cite{csw} studied
Problem \ref{prob-1.1} in the case when $\alpha=0$. The main results in \cite{zjf} are as follows:
 \begin{enumerate}
   \item[$(\mathfrak{i})$](\cite[Lemma 2.3]{zjf})
    If $p\in[1,\infty]$, then $$\|\partial_{\theta}f \|_{\mathcal{H}_{\mathcal{G}}^p(\mathbb{D})}\lesssim \|\dot{F}\|_{\mathcal{L}^{p}(\mathbb{T})};$$
\item[$(\mathfrak{ii})$](\cite[Theorem  1.1]{zjf})
 If $p\in[1,\infty)$, then $$\|\partial_{r}f \|_{\mathcal{B}_{\mathcal{G}}^p(\mathbb{D})}\lesssim \|\dot{F}\|_{\mathcal{L}^{p}(\mathbb{T})};$$
   \item[$(\mathfrak{iii})$](\cite[Theorem  1.2]{zjf})
    If $p\in[1,2)$, then  $$\|\partial_{z}f\|_{  \mathcal{B}_{\mathcal{G}}^p(\mathbb{D})}\lesssim \|\dot{F}\|_{\mathcal{L}^{p}(\mathbb{T})}\;\;\mbox{and}
   \;\;\|\partial_{\overline{z}}f\|_{  \mathcal{B}_{\mathcal{G}}^p(\mathbb{D})}\lesssim \|\dot{F}\|_{\mathcal{L}^{p}(\mathbb{T})};$$
   \item[$(\mathfrak{iv})$](\cite[Theorem 1.3]{zjf})
    If $f$ is a $K$-quasiregular mapping and $p\in[1,\infty]$,
 then $$\|\partial_{z}f\|_{\mathcal{H}_{\mathcal{G}}^p(\mathbb{D})}\lesssim \|\dot{F}\|_{\mathcal{L}^{p}(\mathbb{T})}\;\;\mbox{and}
   \;\;\|\partial_{\overline{z}}f\|_{\mathcal{H}_{\mathcal{G}}^p(\mathbb{D})}\lesssim \|\dot{F}\|_{\mathcal{L}^{p}(\mathbb{T})}.$$
 \end{enumerate}

Here, for two nonnegative quantities $X$ and $Y$,
we denote by   $X\lesssim Y$ if
there exist  two  constants
$M_{1}\geq 1$ and  $M_{2}\geq0 $  such that
$$
X \leq M_{1}Y+M_{2}.
$$

 In \cite{csw}, the second author et al. discussed the generalizations of the arguments in \cite{zjf}. The following are the main results in \cite{csw}:
\begin{enumerate}
  \item[$(\mathfrak{v})$](\cite[Theorem 1.1]{csw})
   If $p\in[1,\infty)$, then $\partial_{z}f$ and $   \partial_{\overline{z}}f \in \mathcal{B}_{\mathcal{G}}^p(\mathbb{D})$;
 if $p=\infty$, then there exists a harmonic mapping $f=\mathcal{K}_{0}[F]$, where $F$ is an absolutely continuous function with $\dot{F}\in \mathcal{L}^{\infty}(\mathbb{T})$, such that neither $\partial_{z}f$ nor $   \partial_{\overline{z}}f $ belongs to $\mathcal{B}_{\mathcal{G} }^{\infty}(\mathbb{D})$;
 \item[$(\mathfrak{vi})$](\cite[Theorem 1.2]{csw})
    If $f$ is a  $(K,K')$-elliptic mapping and $p\in[1,\infty]$,
    then $\partial_{z}f$ and $   \partial_{\overline{z}}f  \in \mathcal{H}_{\mathcal{G}}^p(\mathbb{D})$.
\end{enumerate}
See \cite{geh,    kalaj2019, kalme2011} for more discussions in this line.

\subsection{Main Results}
In this paper, we investigate Problem \ref{prob-1.1} further.
To state our results,
let us introduce the following notations.
\begin{itemize}
  \item Let $\Pi=\{(\alpha, p):\; -1< \alpha< \infty,\; 1\leq p \leq \infty\}$;
  \item Let $\Pi_{1}=\{(\alpha, p):\; 0< \alpha<\infty,\; 1\leq p \leq \infty \}$ $\cup$ $\{(\alpha, p):  \alpha=0,\; 1<p<\infty\}$;
  \item Let $\Pi_{2}=\big\{(\alpha, p):\; -1<\alpha<0,\; 1\leq p<-\frac{1}{\alpha}\big\}$ $\cup$ $ \{ (0,1)\}$;
  \item Let $\Pi_{3}=\{(\alpha, p):\; -1<\alpha<0,\; -\frac{1}{\alpha}\leq p\leq \infty \}$ $\cup$ $\{(0,\infty)\}$.
\end{itemize}
Obviously, $$\Pi=\Pi_{1}\cup\Pi_{2}\cup\Pi_{3}.$$

 First, we determine the range of parameters $\alpha$ and $p$ in which  the
Hardy type  norm
  of  the partial derivatives $\partial_{\theta }f$, $\partial_{r}f$, $\partial_{z}f$ and $ \partial_{\overline{z}}f$ of $f$, i.e.,
$\|\partial_{\theta} f \|_{\mathcal{H}^{p}_{\mathcal{G}}(\mathbb{D})}$, $\|\partial_{r} f \|_{\mathcal{H}^{p}_{\mathcal{G}}(\mathbb{D})}$, $\|\partial_{z} f  \|_{\mathcal{H}^{p}_{\mathcal{G}}(\mathbb{D})}$ and
$\| \partial_{\overline{z}} f  \|_{\mathcal{H}^{p}_{\mathcal{G}}(\mathbb{D})}$,
can be controled by $\|\dot F\|_{ \mathcal{L}^{p}(\mathbb{T})}$.
Our results are as follows.

\bthm\label{thm-1.1}
 Suppose that
 \begin{itemize}
   \item $F$ is an absolutely continuous function on $\mathbb{T}$  and $\dot F\in \mathcal{L}^{p}(\mathbb{T})$ with $p \in [1,\infty]$;
   \item $f=\mathcal{K}_{\alpha}[F]$ in $\mathbb{D}$ with $\alpha\in(-1,\infty)$.
 \end{itemize}
Then the following three statements are true:
\begin{enumerate}
\item\label{thm-1.1-1}
 If $(\alpha, p) \in\Pi$,  then  $\|\partial_{\theta}f\|_{\mathcal{H}_{\mathcal{G}}^p(\mathbb{D})}\lesssim\|\dot{F}\|_{\mathcal{L}^{p}(\mathbb{T})}$;
\item\label{thm-1.1-2}
 If $(\alpha, p) \in\Pi_{1}$, then
 $$
 \|\partial_{r} f  \|_{\mathcal{H}^{p}_{\mathcal{G}}(\mathbb{D})}+ \|\partial_{z} f  \|_{\mathcal{H}^{p}_{\mathcal{G}}(\mathbb{D})}
+ \| \partial_{\overline{z}} f  \|_{\mathcal{H}^{p}_{\mathcal{G}}(\mathbb{D})}
\lesssim \|  F\|_{\mathcal{L}^{\infty}(\mathbb{T})}+ \|\dot F\|_{\mathcal{L}^p(\mathbb{T})};
$$
In particular, if $\alpha=0$ and $p\in(1,\infty)$, then
 $$
 \|\partial_{r} f  \|_{\mathcal{H}^{p}_{\mathcal{G}}(\mathbb{D})}+ \|\partial_{z} f  \|_{\mathcal{H}^{p}_{\mathcal{G}}(\mathbb{D})}
+ \| \partial_{\overline{z}} f  \|_{\mathcal{H}^{p}_{\mathcal{G}}(\mathbb{D})}
\lesssim  \|\dot F\|_{\mathcal{L}^p(\mathbb{T})};
$$
\item\label{thm-1.1-3}
If $(\alpha, p) \in\Pi_{2}\cup \Pi_{3}$ and $f$ is a $(K,K')$-elliptic mapping in $\mathbb{D}$,
then
 $$
\|\partial_{r}f\|_{\mathcal{H}_{\mathcal{G}}^{p}( \mathbb{D} )}
+\|\partial_{z}f\|_{\mathcal{H}_{\mathcal{G}}^{p}( \mathbb{D} )}
+
\|\partial_{\overline{z}}f\|_{\mathcal{H}_{\mathcal{G}}^{p}( \mathbb{D} )}
\lesssim \| \dot F\|_{\mathcal{L}^p(\mathbb{T})}.
$$
\end{enumerate}
\ethm

\bthm\label{thm-1.2}
If $(\alpha, p) \in \Pi_{2}\cup \Pi_{3}$, then there exists a mapping $f $ which satisfies the assumptions of Theorem \ref{thm-1.1}, but
$$
 H_p(\partial_{r}f)=H_p(\partial_{z}f)=H_p(\partial_{\overline{z}}f)=\infty.
$$
\ethm

 The following results are the direct consequences of Theorems  \ref{thm-1.1} and \ref{thm-1.2}.
\bcor\label{cor-1.1}
 Under the assumptions of Theorem \ref{thm-1.1},   the following three statements are true:
\begin{enumerate}
\item If $(\alpha, p) \in\Pi$,  then $\partial_{\theta}f \in\mathcal{H}_{\mathcal{G}}^p(\mathbb{D})$;
\item If $(\alpha, p) \in\Pi_{1}$, then all partial derivatives $\partial_{r}f$, $\partial_{z}f$ and $ \partial_{\overline{z}}f$ of $f$ belong to $\mathcal{H}_{\mathcal{G}}^p(\mathbb{D})$;
\item If $(\alpha, p) \in\Pi_{2}\cup \Pi_{3}$ and $f$ is a $(K,K')$-elliptic mapping in $\mathbb{D}$,
then all partial derivatives $\partial_{r}f$, $\partial_{z}f$ and $ \partial_{\overline{z}}f$ of $f$ belong to $\mathcal{H}_{\mathcal{G}}^p(\mathbb{D})$.
\end{enumerate}
\ecor

\bcor\label{cor-1.2}
 If $(\alpha, p) \in \Pi_{2}\cup \Pi_{3}$, then there exists a mapping $f $ satisfying the assumptions of Theorem \ref{thm-1.1},
but none of the
 partial derivatives $\partial_{r}f$, $\partial_{z}f$ and $ \partial_{\overline{z}}f$ of $f$ belongs to $\mathcal{H}_{\mathcal{G}}^p(\mathbb{D})$.
  \ecor

By Theorem \ref{thm-1.2}, we see that for each pair $(\alpha, p) \in \Pi_{2}\cup \Pi_{3}$,
there exists a mapping $f$ such that all the quantities
 $ H_p(\partial_{r}f)$,
 $ H_p(\partial_{z}f)$ and $H_p(\partial_{\overline{z}}f)$
are infinite.
Notice that
$$B_{p}(\partial_{r}f) \leq H_p(\partial_{r}f),\;\;
 B_{p}(\partial_{z}f) \leq H_p(\partial_{z}f)\;\;
\text{and} \;\;B_{p}(\partial_{\overline{z}}f) \leq  H_p(\partial_{\overline{z}}f).$$
Naturally, one will ask whether the
 partial derivatives $\partial_{r}f$, $\partial_{z}f$ and $ \partial_{\overline{z}}f$ of $f$ belong  to $\mathcal{B}_{\mathcal{G}}^p(\mathbb{D})$
and the corresponding  Bergman type norms can be controlled by $\|\dot F\|_{\mathcal{ L}^{p}(\mathbb{T})}$.
The following results answer  this question.

\bthm\label{thm-1.3}
If $(\alpha, p) \in  \Pi_{2}$ and $f$ satisfies the assumptions of Theorem \ref{thm-1.1},
then
$$
 \|\partial_{r} f  \|_{\mathcal{B}^{p}_{\mathcal{G}}(\mathbb{D})}+ \|\partial_{z} f  \|_{\mathcal{B}^{p}_{\mathcal{G}}(\mathbb{D})}
+ \| \partial_{\overline{z}} f  \|_{\mathcal{B}^{p}_{\mathcal{G}}(\mathbb{D})}
\lesssim   \|  F\|_{\mathcal{L}^{\infty}(\mathbb{T})}+ \|\dot F\|_{\mathcal{L}^p(\mathbb{T})}.
$$
 \ethm

\bthm\label{thm-1.4}
If $(\alpha, p) \in  \Pi_{3}$, then there exists a mapping $f $ which satisfies the assumptions of Theorem \ref{thm-1.1}, but
$$
B_p( \partial_{r}f)=B_p(\partial_{z}f)=B_p(\partial_{\overline{z}}f)=\infty.
$$
\ethm

As the direct consequences of Theorems \ref{thm-1.3} and \ref{thm-1.4}, we have the following.
\bcor\label{cor-1.3}
If $(\alpha, p) \in  \Pi_{2}$ and $f$ satisfies the assumptions of Theorem \ref{thm-1.1},
then all partial derivatives $\partial_{r}f$, $\partial_{z}f$ and $\partial_{\overline{z}}f$ belong  to $\mathcal{B}_{\mathcal{G}}^{p}(\mathbb{D})$.
\ecor

\bcor\label{cor-1.4}
If $(\alpha, p) \in  \Pi_{3}$, then there exists a mapping $f $ satisfying the assumptions of Theorem \ref{thm-1.1},
but none of
 the partial derivatives $\partial_{r}f$, $\partial_{z}f$ and $\partial_{\overline{z}}f$ of $f$ belongs to $\mathcal{B}_{\mathcal{G}}^{p}(\mathbb{D})$.
\ecor

\begin{rem}\label{remark-1.1}
The three examples constructed in Section \ref{sec-4}
show that the assumption that $f$ is a $(K,K')$-elliptic mapping
in Theorem \ref{thm-1.1}(3) and in Corollary \ref{cor-1.1}(3) can not be removed.
\end{rem}

\br\label{remark-1.2} Even for the special case when $\alpha=0$, we generalize the discussions in both \cite{csw} and \cite{zjf}. The details are as follows.
\begin{itemize}
\item[$(\mathfrak{a})$]
When $\alpha=0$ and $p\in[1,\infty]$, Theorem  \ref{cor-1.1}(1) reduces to \cite[Lemma 2.3]{zjf}.
\item[$(\mathfrak{b})$]
When $\alpha=0$ and $p\in(1,\infty)$,
Theorem \ref{thm-1.1}(2) (resp. Corollary \ref{cor-1.1}(2)) is an improvement of \cite[Theorems 1.1$-$1.3]{zjf} (resp. \cite[Theorems 1.1(1) and 1.2]{csw}). Also,
Theorem \ref{thm-1.1}(2) and Corollary \ref{cor-1.1}(2) show that the assumption that $f$ is a $(K,K')$-elliptic mapping in \cite[Theorem  1.3]{zjf} and \cite[Theorem 1.2]{csw} can be removed.
\item[$(\mathfrak{c})$]
When $\alpha=0$ and $p=\infty$, Corollary \ref{cor-1.4} reduces to \cite[Theorem 1.1(2)]{csw}.
\item[$(\mathfrak{d})$]
When $\alpha=0$ and $p=1$, Theorem \ref{thm-1.3} (resp. Corollary \ref{cor-1.3}) reduces to \cite[Theorems  1.1 and 1.2]{zjf} (resp. \cite[Theorem 1.1(1)]{csw}).
\item[$(\mathfrak{e})$]
When $\alpha=0$ and $p\in\{1,\infty\}$, Corollary \ref{cor-1.1}(3) reduces to \cite[Theorem 1.2]{csw}.
\end{itemize}
\er

This paper is organized as follows.
 In Section 2, necessary terminologies will be introduced, and several known results will be recalled.
 In Section 3, first, a series of lemmas will be established. Based on these lemmas, Theorems \ref{thm-1.1} and \ref{thm-1.3} will be proved, and in Section \ref{sec-4}, first, three examples will be constructed, and then,  Theorems \ref{thm-1.2} and \ref{thm-1.4} will be shown.


\section{Preliminaries }\label{sec-2}
In this section, we shall recall some necessary terminology and useful known results.
We start with the definition of convex functions.

Let $I$ be an interval in $\mathbb{R}$. A function $f:I\rightarrow\mathbb{R}$ is said to be \emph{convex} if for all $x$, $y\in I$ and $\lambda\in[0,1]$,
$$
f\big(\lambda x+(1-\lambda)y\big)\leq\lambda f(x)+(1-\lambda)f(y).
$$

\subsection*{Jensen's inequality.} Suppose that $\varphi:[\alpha,\beta]\rightarrow\mathbb{R}$   is a convex function,
$f$ and $p$ are  integrable in $[a,b]$, where $\alpha<\beta$ and $a<b$.
If for any $x\in [a,b]$, $f(x)\in[\alpha,\beta]$, $p(x)\geq 0$ and
$\int_{a}^{b}p(x)dx>0$,
then
$$\varphi\left(\frac{\int_{a}^{b}f(x)p(x)\,dx}{\int_{a}^{b}p(x)dx}\right)
\leq \frac{\int_{a}^{b}\varphi\big(f(x)\big)p(x)dx}{\int_{a}^{b}p(x)dx}.$$

\subsection*{Gauss hypergeometric functions.}
For any $a\in\mathbb{R}$ and $k\in \mathbb{N}_{0}$,
let
$$
(a)_{k}=\left\{\begin{array}{ll}
1, & \text{if} \;\;\;k=0,\\
 a(a+1)\ldots (a+k-1), & \text{if} \;\;\;k\geq 1,
 \end{array}\right.
 $$ which is called the {\it factorial function}.

If $a$ is neither zero nor a negative integer, then
$$
 (a)_{k}=\frac{\Gamma(a+k)}{\Gamma(a)}
$$
(cf. \cite[p.23]{rain}).
For  $x\in \mathbb{R}$,
the {\it Gauss hypergeometric function} or the {\it hypergeometric series} is defined as follows:
\be\label{eq-2.1}
 _{2}F_{1}(a,b;c;x)
=\sum_{k=0}^{\infty}\frac{(a)_{k}(b)_{k}}{k!(c)_{k}}x^{k},
\ee
where $a,b\in \mathbb{R}$  and $c$ is neither zero nor a negative integer (cf. \cite[p.46]{rain}).

If $c-a-b>0$, then the series $_{2}F_{1}(a,b;c;x)$ is absolutely convergent for all $x$ with $|x|\leq 1$ (cf. \cite[Section 31]{rain}), and
\begin{eqnarray}\label{eq-2.2}
 _{2}F_{1}\left(a,b;c;1\right)
=\frac{\Gamma(c)\Gamma(c-a-b)}{\Gamma(c-a)\Gamma(c-b)}
\end{eqnarray}
(cf. \cite[p.49]{rain}).

If $c-a-b<0$ and $|x|<1$, then
\begin{eqnarray}\label{eq-2.3}
 _{2}F_{1}(a ,b ;c ;x)=(1-x)^{c-a-b}\;_{2}F_{1}(c-a ,c-b ;c ;x)
\end{eqnarray}(cf.
\cite[p.60]{rain}).

Elementary computations guarantee the following useful formula:
\begin{eqnarray*}
\frac{d}{dx}\;_{2}F_{1}(a,b;c;x)
= \frac{ab}{c}
 \;_{2}F_{1}(a+1,b+1;c+1;x).
\end{eqnarray*}

\subsection*{Distributions}
A distribution $F\in \mathcal{D}'^{m}(\mathbb{T})$ on $\mathbb{T}$ of order less than or equal to $m\in \mathbb{N}_{0}$ is a linear
form on $\mathcal{C}^{\infty}(\mathbb{T})$ such that there exists a constant $C>0$ such that for any $\varphi\in \mathcal{C}^{\infty}(\mathbb{T})$,
$$
|\langle F,\varphi\rangle|\leq C\|\varphi\|_{\mathcal{C}^{m}(\mathbb{T})},
$$
where $\|\cdot\|_{\mathcal{C}^{m}(\mathbb{T})}$ denotes the norm of $m$-times continuously differentiable functions on $\mathbb{T}$,
$\langle\cdot ,\cdot \rangle$ denotes the distributional pairing
 and
\begin{eqnarray}\label{eq-2.4}
\langle F,\varphi\rangle:=\frac{1}{2\pi}\int_{\mathbb{T}} \varphi(e^{i\theta})F(e^{i\theta})d\theta
\end{eqnarray}(cf. \cite[p.570]{tade}).
Then a  function $F\in \mathcal{L}^1(\mathbb{T})$ is identified with the distribution of order 0 (cf. \cite[Section 3]{aj}).
Notice that the space $ \mathcal{D}'^{m}(\mathbb{T})$ is naturally identified with the dual of $\mathcal{C}^{m}(\mathbb{T})$.

A distribution $F\in \mathcal{D}'(\mathbb{T})$ on $\mathbb{T}$ is an element in $\mathcal{D}'^{m}(\mathbb{T})$  for some $m\in \mathbb{N}_0$, that is,
$$ \mathcal{D}'(\mathbb{T})=\cup_{m\geq 0}\mathcal{D}'^{m}(\mathbb{T})$$ (cf. \cite[Section 3]{aj}).
Obviously,
\begin{eqnarray*}
\mathcal{L}^{\infty}(\mathbb{T})\subset \mathcal{L}^{1}(\mathbb{T}) =\mathcal{D}'^{0}(\mathbb{T}) \subset\mathcal{D}'(\mathbb{T}).
\end{eqnarray*}
The space $\mathcal{D}'(\mathbb{T})$ of distributions is topologised by means
of the semi-norms
$$\mathcal{D}' (\mathbb{T})\ni F\mapsto \big|\langle F ,\varphi\rangle\big|$$
 for every  $\varphi\in \mathcal{C}^{\infty}(\mathbb{T})$.
We say that $F_{k}\rightarrow F$ in $\mathcal{D}' (\mathbb{T})$ as $k\rightarrow\infty$
if
$\langle F_{k},\varphi\rangle\rightarrow\langle F ,\varphi\rangle $ as $k\rightarrow\infty$ for every  $\varphi\in \mathcal{C}^{\infty}(\mathbb{T})$
 (cf.  \cite[Section 3]{aj}).

Next, let us recall several known results, which will be used
later on.

\begin{Lem}\label{Lem-H}{\rm (cf. \cite[p.158]{kuang})}
Let $a_{1},\ldots,a_{n}\in \mathbb{R}$, where $n\in \mathbb{N}$.
Then
$$
\left(\sum_{k=1}^{n}|a_{k}|\right)^{p}
\leqslant\left\{\begin{array}{ll}
\sum_{k=1}^{n}|a_{k}|^{p} , & \text{if} \;\;\;0<p<1,\\
 n^{p-1}\sum_{k=1}^{n}|a_{k}|^{p}, & \text{if} \;\;\; p \geqslant 1.
 \end{array}\right.
 $$
\end{Lem}

\begin{Lem}\label{Lem-I}{\rm (\cite[Lemma $2.3$]{aj})}
Let $\alpha>0$ and $r\in [0,1)$. Then
$$
\frac{1}{2\pi}\int_{0}^{2\pi} \frac{(1-r^2)^\alpha}{| 1-re^{i\theta}| ^{\alpha+1}}\,d\theta
\leq\frac{\Gamma(\alpha)}{\Gamma^2(\frac{\alpha+1}{2})}.
$$
\end{Lem}


\begin{Lem}\label{Lem-J}{\rm (cf. \cite[p.18-19]{rain})}
For $s$, $t\in(-1,\infty)$, we have
 $$
\int_{0}^{1} (1-r )^{s}r^{t}\,dr=\frac{\Gamma(s+1) \Gamma(t+1)}{\Gamma(s+t+2)}.
$$
\end{Lem}

\section{Proofs  of Theorems \ref{thm-1.1} and \ref{thm-1.3} }\label{sec-3}
 The aim of this section is to prove Theorems \ref{thm-1.1} and \ref{thm-1.3}. The proofs are based on the following lemmas.

\subsection{Several lemmas}

 \begin{lem}\label{lem-3.1}
 Suppose that
 \begin{itemize}
   \item $F$ is an absolutely continuous function on $\mathbb{T}$ and $\dot F\in \mathcal{L}^{p}(\mathbb{T})$ with $p \in [1,\infty]$;
   \item $f=\mathcal{K}_{\alpha}[F]$ in $\mathbb{D}$ with $\alpha\in(-1,\infty)$.
 \end{itemize}
Then the following two statements are true:
\begin{enumerate}
  \item  For all $\alpha\in(-1,\infty)$ and $p \in [1,\infty]$,
  $$\|\partial_{\theta}f\|_{\mathcal{H}_{\mathcal{G}}^p(\mathbb{D})}\lesssim\|\dot{F}\|_{\mathcal{L}^{p}(\mathbb{T})};$$
  \item  For $\alpha=0$ and all $p\in(1,\infty)$,
  $$\|  \partial_{r}f\|_{\mathcal{H}_{\mathcal{G}}^p(\mathbb{D})}
  +\|  \partial_{\overline{z}}f\|_{\mathcal{H}_{\mathcal{G}}^p(\mathbb{D})}
  +\|  \partial_{\overline{z}}f\|_{\mathcal{H}_{\mathcal{G}}^p(\mathbb{D})}
 \lesssim\|\dot F \|_{\mathcal{L}^{p}(\mathbb{T})}.$$
\end{enumerate}

\end{lem}
\bpf (1) In order to prove the first statement in the lemma, we consider the case when $p\in[1,\infty)$ and the case when $p=\infty$,
separately.
\bca\label{case-1}  Suppose that $p\in[1,\infty)$.
\eca
For $z=re^{i\theta}\in \mathbb{D}$  and $\alpha\in(-1,\infty)$,
since the  mappings
$$(z,t)\mapsto \mathcal{K}_{\alpha}(ze^{-it})F(e^{it})
\;\text{and}\;
(z,t)\mapsto  \frac{\partial  \mathcal{K}_{\alpha}(ze^{-it})}{\partial \theta}F(e^{it})
 $$
are continuous in $\overline{\mathbb{D}}_{\rho_{0}}\times[0,2\pi]$,
where $\rho_{0}\in(0,1)$ and
$\overline{\mathbb{D}}_{\rho_{0}}=\{z\in \mathbb{C}:|z|\leq \rho_{0}\}$,
 we obtain from \eqref{eq-1.3} that
\begin{align*}
\begin{split}
\partial_{\theta}f(re^{i\theta})
&= \frac{1}{2\pi} \int_{0}^{2\pi}  \Big( \frac{\partial}{\partial \theta}  \mathcal{K}_{\alpha}(re^{i(\theta-t)}) \Big)F(e^{it}) \,dt
 = -\frac{1}{2\pi} \int_{0}^{2\pi}  \Big( \frac{\partial}{\partial t}   \mathcal{K}_{\alpha}(re^{i(\theta-t)})  \Big)F(e^{it})\,dt.
\end{split}
\end{align*}
Recall that $F$ is absolutely continuous on $\mathbb{T}$ with $\dot{F}\in \mathcal{L}^{p}(\mathbb{T})$.
Then the integration by parts
  leads to
\begin{align}\label{eq-3.1}
\begin{split}
\partial_{\theta}f(re^{i\theta})
= \frac{1}{2\pi}\int_{0}^{2\pi}\mathcal{K}_{\alpha}(re^{i(\theta-t)})\dot F(e^{it})\,dt.
\end{split}
\end{align}
For $\alpha\in(-1,\infty)$,  let
\begin{equation}\label{eq-3.2}
I_{\alpha,1}(r):=\int_{0}^{2\pi}\frac{(1-r^2)^{\alpha+1}}{|  1-re^{i t}|^{\alpha+2}}\,dt,
\end{equation} where $r\in [0,1)$.
By Lemma I, we get
\begin{equation}\label{eq-3.3}
I_{\alpha,1}(r)\leq \frac{2\pi\Gamma(\alpha+1)}{ \Gamma^2(\frac{\alpha}{2}+1) }.
\end{equation}
By \eqref{eq-1.4} and \eqref{eq-3.2}, we easily find that
\begin{align}\label{eq-3.4}
\begin{split}
 \int_{0}^{2\pi} \mathcal{K}_{\alpha}(re^{i(\theta-t)})\,d\theta
 = \mathcal{C}_{\alpha} I_{\alpha,1}(r),
\end{split}
\end{align}
where $\mathcal{C}_{\alpha}$ is the constant from \eqref{eq-1.5}.
Then \eqref{eq-3.1} and Jensen's inequality ensure  that
\begin{align}\label{eq-3.5}
\begin{split}
| \partial_{\theta}f(re^{i\theta})|^p
&\leq\frac{\big(\mathcal{C}_{\alpha} I_{\alpha,1} (r)\big)^{p}}{(2\pi)^p}\left(\int_{0}^{2\pi}\frac{\mathcal{K}_{\alpha}(re^{i(\theta-t)})}{\mathcal{C}_{\alpha}I_{\alpha,1}(r)}|\dot F(e^{it})|\,dt\right)^p\\
&\leq\frac{\big(\mathcal{C}_{\alpha} I_{\alpha,1} (r)\big)^{p-1}}{(2\pi)^p}\int_{0}^{2\pi}\mathcal{K}_{\alpha}(re^{i(\theta-t)})|\dot F(e^{it})|^p\,dt.
\end{split}
\end{align}

For fixed $r\in [0,1)$, since the assumption of $\dot F\in \mathcal{L}^{p}(\mathbb{T})$ ensures that
$$
\int_{0}^{2\pi}\int_{0}^{2\pi}  \mathcal{K}_{\alpha}(re^{i(\theta-t)})|\dot F(e^{it})|^p    d\theta dt
\leq\frac{2^{\alpha+3}\pi^{2} \mathcal{C}_{\alpha}}{ 1-r } \| \dot{F} \|^{p}_{\mathcal{L}^{p}(\mathbb{T})}
<\infty,
$$
by \eqref{eq-3.5} and  Fubini's Theorem, we obtain that
\begin{align}\label{eq-3.6}
\begin{split}
\int_{0}^{2\pi}| \partial_{\theta}f(re^{i\theta})| ^{p} \,d\theta
\leq&\frac{\big(\mathcal{C}_{\alpha}I_{\alpha,1} (r)\big )^{p-1}}{(2\pi)^p}\int_{0}^{2\pi}\,d\theta \int_{0}^{2\pi} \mathcal{K}_{\alpha}(re^{i(\theta-t)})|\dot F(e^{it})|^p\,dt\\
=&\frac{\big(\mathcal{C}_{\alpha}I_{\alpha,1} (r)\big )^{p-1} }{(2\pi)^p}\int_{0}^{2\pi}|\dot F(e^{it})|^p\,dt \int_{0}^{2\pi} \mathcal{\mathcal{K}}_{\alpha}(re^{i(\theta-t)})\,d\theta.
\end{split}
\end{align}
Therefore,
\begin{align}\label{eq-3.7}
\begin{split}
\quad\int_{0}^{2\pi}| \partial_{\theta}f(re^{i\theta})| ^{p} \,d\theta
\leq&\frac{\big( \mathcal{C}_{\alpha}I_{\alpha,1}(r)\big)^p}{(2\pi)^{p-1 }}\|\dot F\|_{\mathcal{L}^p(\mathbb{T})}^{p}
\qquad\;\quad\;\;\;\;\text{(by \eqref{eq-3.4}  and \eqref{eq-3.6})}\\
\leq& 2\pi \|\dot F\|_{\mathcal{L}^p(\mathbb{T})}^{p},\qquad\qquad\qquad\qquad\;\text{(by \eqref{eq-1.5} and \eqref{eq-3.3})}
\end{split}
\end{align}
which implies that for any $\alpha\in(-1,\infty)$ and $p\in[1,\infty)$,
\begin{equation*}
 \| \partial_{\theta}f\|_{\mathcal{H}_{\mathcal{G}}^p(\mathbb{D})}=\sup\limits_{r\in [0,1)} \left(\frac{1}{2\pi} \int_{0}^{2\pi}| \partial_{\theta}f(re^{i\theta})|^p\,d\theta \right)^{\frac{1}{p}}\leq \|\dot F \|_{\mathcal{L}^{p}(\mathbb{T})} <\infty.
\end{equation*}

\smallskip
\bca\label{case-2} Suppose that $p=\infty$.
\eca
For $z=re^{i\theta}\in \mathbb{D}$ and $\alpha\in(-1,\infty)$, it follows from \eqref{eq-1.4}, \eqref{eq-1.5} and Lemma I
that
$$
\frac{1}{2\pi}\int_{0}^{2\pi} \mathcal{K}_{\alpha}(re^{it})dt\leq 1.
$$
Since $\dot F\in \mathcal{L}^{\infty}(\mathbb{T})$, we infer from \eqref{eq-3.1} and the above inequality that
\begin{align*}
\begin{split}
\big|\partial_{\theta}f(re^{i\theta})\big|
\leq\frac{1}{2\pi}\int_{0}^{2\pi}\mathcal{K}_{\alpha}(re^{i(\theta-t)})\big|\dot F(e^{it})\big|\,dt
\leq \|\dot F\|_{\mathcal{L}^{\infty}(\mathbb{T})},
\end{split}
\end{align*}
which gives $$ \| \partial_{\theta}f\|_{\mathcal{H}_{\mathcal{G}}^{\infty}(\mathbb{D})} \leq \|\dot F \|_{\mathcal{L}^{\infty}(\mathbb{T})}.$$
Thus, the first statement in the lemma is true.

(2) If $\alpha=0$ and $p\in(1,\infty)$, then $f=\mathcal{K}_{\alpha}[F]$ is harmonic in $\mathbb{D}$ (cf. \cite[p.41]{pav}),
and so,
$f$ admits a decomposition $f=h+\overline{g}$, where $h$ and $g$ are analytic  in $\mathbb{D}$.
For $z=re^{i\theta}\in \mathbb{D}$, obviously,
\begin{align}\label{eq-3.8}
\begin{split}
r \partial_{r}f(z)
=zh'(z)+\overline{zg'(z)}
 \;\;\text{and}\;\;
 \partial_{\theta}f(z)=izh'(z)+\overline{izg'(z)}.
\end{split}
\end{align}
Let
$$
izh'(z)=u_{1}(z) +iv_{1}(z)\;\;\text{and}\;\;izg'(z)=u_{2}(z) +iv_{2}(z),
$$
where $u_{j}$ and $v_{j}$ are real harmonic mappings in $\mathbb{D}$
for $j\in\{1,2\}$.
Obviously,
$u_{1}(0)=v_{1}(0)=u_{2}(0)=v_{2}(0)=0$.
Let
$$
\Phi (z)=izh'(z) + izg'(z)
$$
  in $\mathbb{D}$.
Then
\begin{align}\label{eq-3.9}
\begin{split}
\partial_{\theta}f(z)= \big(u_{1}(z) +u_{2} (z)\big  )+i \big(v_{1}(z) -v_{2}(z)\big )
\end{split}
\end{align}
and
$$
\Phi(z)= \big  (u_{1}(z) +u_{2}(z)  \big  )+i \big  (v_{1}(z) +v_{2}(z)\big   ).
$$

Since Lemma \ref{lem-3.1}(1) ensures that
 $$ \| \partial_{\theta}f\|_{\mathcal{H}_{\mathcal{G}}^p(\mathbb{D})}\lesssim\|\dot F \|_{\mathcal{L}^{p}(\mathbb{T})},$$
we know from \eqref{eq-3.9}  that
\begin{align}\label{eq-3.10}
\begin{split}
  \| u_{1}+u_{2}\|_{\mathcal{H}_{\mathcal{G}}^p(\mathbb{D})} +\| v_{1} -v_{2}\|_{\mathcal{H}_{\mathcal{G}}^p(\mathbb{D})}
 \leq 2\| \partial_{\theta}f\|_{\mathcal{H}_{\mathcal{G}}^p(\mathbb{D})}
 \lesssim\|\dot F \|_{\mathcal{L}^{p}(\mathbb{T})}.
\end{split}
\end{align}
 Notice that $\Phi $ is analytic in $\mathbb{D}$. Then
$v_{1} +v_{2} $ is a harmonic conjugate of $u_{1} +u_{2}$ with $v_{1} (0) +v_{2} (0)=0$.
Since $ \| u_{1}+u_{2}\|_{\mathcal{H}_{\mathcal{G}}^p(\mathbb{D})}
 \lesssim\|\dot F \|_{\mathcal{L}^{p}(\mathbb{T})}$,
 we obtain from   Theorem E that
\begin{align}\label{eq-3.11}
\begin{split}
\| v_{1}+v_{2}\|_{\mathcal{H}_{\mathcal{G}}^p(\mathbb{D})}
\lesssim \| u_{1}+u_{2}\|_{\mathcal{H}_{\mathcal{G}}^p(\mathbb{D})}
 \lesssim\|\dot F \|_{\mathcal{L}^{p}(\mathbb{T})},
\end{split}
\end{align}
which, together with  Lemma H or Minkowski's inequality, implies that
 $$
 \|\Phi \|_{\mathcal{H}_{\mathcal{G}}^p(\mathbb{D})}
\lesssim \| u_{1}+u_{2}\|_{\mathcal{H}_{\mathcal{G}}^p(\mathbb{D})}+\| v_{1}+v_{2}\|_{\mathcal{H}_{\mathcal{G}}^p(\mathbb{D})}
 \lesssim\|\dot F \|_{\mathcal{L}^{p}(\mathbb{T})}.
 $$
 Notice that \eqref{eq-3.10} and \eqref{eq-3.11} ensure that
 $$\| v_{1}-v_{2}\|_{\mathcal{H}_{\mathcal{G}}^p(\mathbb{D})}+\| v_{1}+v_{2}\|_{\mathcal{H}_{\mathcal{G}}^p(\mathbb{D})}
 \lesssim\|\dot F \|_{\mathcal{L}^{p}(\mathbb{T})}.$$
Again, Lemma H or Minkowski's inequality implies that
\begin{align}\label{eq-3.12}
\begin{split}
\| v_{1}\|_{\mathcal{H}_{\mathcal{G}}^p(\mathbb{D})}+\| v_{2}\|_{\mathcal{H}_{\mathcal{G}}^p(\mathbb{D})}
 \lesssim\|\dot F \|_{\mathcal{L}^{p}(\mathbb{T})}.
\end{split}
\end{align}

Further, since $ zg'(z)=v_{2}(z) -iu_{2}(z) $, we infer from  Theorem E, Lemma H or Minkowski's inequality and \eqref{eq-3.12} that
\begin{align*}
\begin{split}
\| zg'\|_{\mathcal{H}_{\mathcal{G}}^p(\mathbb{D})}
 \lesssim\| u_{2}\|_{\mathcal{H}_{\mathcal{G}}^p(\mathbb{D})}+\| v_{2}\|_{\mathcal{H}_{\mathcal{G}}^p(\mathbb{D})}
  \lesssim \| v_{2}\|_{\mathcal{H}_{\mathcal{G}}^p(\mathbb{D})}
 \lesssim\|\dot F \|_{\mathcal{L}^{p}(\mathbb{T})}.
\end{split}
\end{align*}
Observe that $g'$ is continuous in $\overline{\mathbb{D}}_{\frac{1}{2}}$ and $|g'(z)|\leq 2|zg'(z)|$
in $\mathbb{D}\backslash\overline{\mathbb{D}}_{\frac{1}{2}}$.
Then we have that
 $$\| g'\|_{\mathcal{H}_{\mathcal{G}}^p(\mathbb{D})}
 \lesssim\|\dot F \|_{\mathcal{L}^{p}(\mathbb{T})}.$$

By \eqref{eq-3.8}, the fact $$\| zg'\|_{\mathcal{H}_{\mathcal{G}}^p(\mathbb{D})}+\| \partial_{\theta}f\|_{\mathcal{H}_{\mathcal{G}}^p(\mathbb{D})}
 \lesssim\|\dot F \|_{\mathcal{L}^{p}(\mathbb{T})} $$ and Lemma  H or Minkowski's  inequality,
we find that
$$ \| zh'\|_{\mathcal{H}_{\mathcal{G}}^p(\mathbb{D})}
 \lesssim\|\dot F \|_{\mathcal{L}^{p}(\mathbb{T})} .$$
Since $h'$ is continuous in $\overline{\mathbb{D}}_{\frac{1}{2}}$ and $|h'(z)|\leq 2|zh'(z)|$
in $\mathbb{D}\backslash\overline{\mathbb{D}}_{\frac{1}{2}}$,
 we have that
 $$\| h'\|_{\mathcal{H}_{\mathcal{G}}^p(\mathbb{D})}
 \lesssim\|\dot F \|_{\mathcal{L}^{p}(\mathbb{T})}.$$

Then \eqref{eq-3.8}, the fact
$$\|  h'\|_{\mathcal{H}_{\mathcal{G}}^p(\mathbb{D})}+\|  g'\|_{\mathcal{H}_{\mathcal{G}}^p(\mathbb{D})}
 \lesssim\|\dot F \|_{\mathcal{L}^{p}(\mathbb{T})} $$ and  Lemma  H  or Minkowski's  inequality
  imply that
 $$\| \partial_{r}f\|_{\mathcal{H}_{\mathcal{G}}^p(\mathbb{D})}
 \lesssim\|\dot F \|_{\mathcal{L}^{p}(\mathbb{T})} .$$
The proof of the lemma is complete.
\epf

The following is an auxiliary result.
\begin{lem}\label{lem-3.2}
For $\alpha\in(-1,0)$ and $r\in[\frac{1}{2},1)$, we have
 $$
 \int_{0}^{2\pi}\frac{dt}{|1-re^{it}|^{\alpha+1}}
 \leq  \frac{3^{\frac{\alpha+1}{2}} }{2^{ \alpha-1 } }  \frac{\Gamma(-\alpha) \Gamma(\frac{1}{2})}{\Gamma(\frac{1}{2}-\alpha)}.
 $$
\end{lem}
\bpf
For any $r\in[\frac{1}{2},1)$ and $t\in[0,2\pi]$,
since
$$
\frac{1+3r^{2}}{6r-2}
=\frac{1}{2}\left(  r-\frac{1}{3} +\frac{4}{9\left(r-\frac{1}{3}\right) } \right)+\frac{1}{3}
>1\geq\cos t,
$$
we see that
$$ 1+r^{2}-2r\cos t\geq \frac{2 }{3}(1-\cos t),$$
and so,

\begin{align}\label{eq-3.13}
\begin{split}
 \int_{0}^{2\pi}\frac{dt}{|1-re^{it}|^{\alpha+1}}
&=2\int_{0}^{ \pi}\frac{dt}{(1+r^{2}-2r\cos
t)^{\frac{\alpha+1}{2}}} \leq  \frac{3^{\frac{\alpha+1}{2}}
}{2^{\frac{\alpha-1}{2}} }   \int_{0}^{ \pi}\frac{dt}{(1-\cos
t)^{\frac{\alpha+1}{2}}}\\
&=\frac{3^{\frac{\alpha+1}{2}}
}{2^{ \alpha-1 } }  \int_{0}^{ \frac{\pi}{2}}\frac{dt}{(\sin
t)^{\alpha+1}}.
\end{split}
\end{align}



For $\alpha\in(-1,0)$, it follows from Lemma J  that
\be\label{eq-3.14}
 \int_{0}^{ \frac{\pi}{2}}\frac{dt}{(\sin t)^{\alpha+1}}
 =\int_{0}^{1} \frac{ds}{s^{\alpha+1} \sqrt{1-s^{2}}}
\leq \int_{0}^{1} \frac{ds}{s^{\alpha+1} \sqrt{1-s }}
=\frac{\Gamma(-\alpha)
\Gamma(\frac{1}{2})}{\Gamma(\frac{1}{2}-\alpha)}.
 \ee
Then the combination of (\ref{eq-3.13}) and (\ref{eq-3.14}) gives
\begin{align*}
\begin{split}
 \int_{0}^{2\pi}\frac{dt}{|1-re^{it}|^{\alpha+1}}
\leq\frac{3^{\frac{\alpha+1}{2}} }{2^{ \alpha-1 } }  \int_{0}^{ \frac{\pi}{2}}\frac{dt}{(\sin t)^{\alpha+1}}
\leq \frac{3^{\frac{\alpha+1}{2}} }{2^{ \alpha-1 } }  \frac{\Gamma(-\alpha) \Gamma(\frac{1}{2})}{\Gamma(\frac{1}{2}-\alpha)},
\end{split}
\end{align*}
which is what we need.
\epf

 \medskip

Assume that $F\in \mathcal{L}^{\infty}(\mathbb{T})$.
For $\alpha\in(-1,\infty)$, let
\begin{equation}\label{eq-3.15}
  \mathscr{J}_{\alpha,1}(z):=\frac{\alpha \mathcal{C}_{\alpha}}{2\pi}\int_{0}^{2\pi}\frac{(1-| z |^2)^{\alpha+1}}{| e^{it}-z |^{\alpha+2}}F(e^{it})\,dt
\end{equation}
and
\begin{align}\label{eq-3.16}
\begin{split}
  \mathscr{J}_{\alpha,2}(z)
  :=&\frac{i\mathcal{C}_{\alpha}(1-| z |^2)^\alpha}{2\pi}\int_{0}^{2\pi} e^{-  \frac{i \alpha   t}{2}}  \Big(  \frac{\partial}{\partial{t}}
\frac{1}{(e^{-it}-\overline{z})^{\frac{\alpha}{2}}(1-ze^{-it})^{\frac{\alpha}{2}+1}} \Big)F(e^{it})\,dt\\
&-\frac{i\mathcal{C}_{\alpha}(1-| z |^2)^\alpha}{2\pi}\int_{0}^{2\pi}   e^{  \frac{i \alpha   t}{2}}   \Big(     \frac{\partial}{\partial{t}}
\frac{1}{(e^{it}-z)^{\frac{\alpha}{2}}(1-\overline{z}e^{it})^{\frac{\alpha}{2}+1}}  \Big)  F(e^{it})\,dt
\end{split}
\end{align} in $\mathbb{D}$,
where $\mathcal{C}_{\alpha}$ is the constant from \eqref{eq-1.5}.

Based on Lemma \ref{lem-3.2}, we have the following estimates concerning $ \mathscr{J}_{\alpha,1}$ and $\mathscr{J}_{\alpha,2}$.
\begin{lem}\label{lem-3.3}
Under the assumptions of Lemma \ref{lem-3.1}, the following three statements are true:
\begin{enumerate}
\item For all $\alpha\in(-1,\infty)$ and $p\in [1,\infty)$
or for all $\alpha\in(0,\infty)$ and $p=\infty$,
$$
\|\mathscr{J}_{\alpha,1}\|_{\mathcal{H}_{\mathcal{G}}^{p}(\mathbb{D})}\lesssim\|F\|_{\mathcal{L}^{\infty}(\mathbb{T})};
$$
\item   For all $\alpha\in(0,\infty)$ and $p\in[1,\infty]$,
$$
\|\mathscr{J}_{\alpha,2}\|_{\mathcal{H}_{\mathcal{G}}^{p}(\mathbb{D})}
\lesssim
\|F\|_{\mathcal{L}^{\infty}(\mathbb{T})}+\|\dot{F}\|_{\mathcal{L}^{p}(\mathbb{T})};
$$
\item  For $\alpha=0$ and all $p\in[1,\infty)$ or for all $\alpha\in(-1,0)$ and
$p\in\big[1,-\frac{1}{\alpha}\big)$,
$$
\|\mathscr{J}_{\alpha,2}\|_{\mathcal{B}_{\mathcal{G}}^{p}(\mathbb{D})}\lesssim\|F\|_{\mathcal{L}^{\infty}(\mathbb{T})}+\|\dot{F}\|_{\mathcal{L}^{p}(\mathbb{T})}.
$$
\end{enumerate}

\end{lem}
\bpf
In the following, we divide the discussions into two cases.

\bca\label{case-3}  Suppose that $p=\infty$.
\eca
First, we estimate the quantity
$\|\mathscr{J}_{\alpha,1}\|_{\mathcal{H}_{\mathcal{G}}^{\infty}(\mathbb{D})}$ for $\alpha\in(0,\infty)$.
Since $ \|F\|_{\mathcal{L}^{\infty}(\mathbb{T})}<\infty$, we deduce from \eqref{eq-1.5}, (\ref{eq-3.15}) and Lemma I that for $z\in \mathbb{D}$,
\begin{align*}
\begin{split}
|\mathscr{J}_{\alpha,1}(z)|
\leq \frac{\alpha\mathcal{C}_{\alpha}  \| F\|_{\mathcal{L}^{\infty}(\mathbb{T})} }{2\pi}
\int_{0}^{2\pi}\frac{(1-|z|^{2})^{\alpha+1}}{ |1-re^{it}|^{\alpha+2}}dt
\leq  \alpha \|F\|_{\mathcal{L}^{\infty}(\mathbb{T})},
\end{split}
\end{align*}
which shows that
\begin{align}\label{eq-3.17}
\begin{split}
\|\mathscr{J}_{\alpha,1}\|_{\mathcal{H}_{\mathcal{G}}^{\infty}(\mathbb{D})}\lesssim\|F\|_{\mathcal{L}^{\infty}(\mathbb{T})}.
\end{split}
\end{align}

Next, we estimate the quantity
$\|\mathscr{J}_{\alpha,2}\|_{\mathcal{H}_{\mathcal{G}}^{\infty}(\mathbb{D})}$ for $\alpha\in(0,\infty)$.
In fact, for $\alpha\in(-1,\infty)$ and $z=re^{i\theta}\in \mathbb{D}$, since $F$ is continuous,
 we deduce from \eqref{eq-3.16} and the integration by parts that
\begin{align*}
\begin{split}
 \mathscr{J}_{\alpha,2}(z)
=&\frac{\mathcal{C}_{\alpha}(1-| z |^2)^\alpha}{ 2\pi}\int_{0}^{2\pi}
  \frac{\frac{\mathrm{d}}{\mathrm{d}t}\big(i e^{\frac{i\alpha t}{2}}  F(e^{it})\big)}
{(e^{it}-z)^{\frac{\alpha}{2}}(1-\overline{z}e^{it})^{\frac{\alpha}{2}+1}} \,dt\\
&-\frac{\mathcal{C}_{\alpha}(1-| z |^2)^\alpha}{ 2\pi}\int_{0}^{2\pi}
 \frac{ \frac{\mathrm{d}}{\mathrm{d}t}\big( ie^{-\frac{i \alpha   t}{2}}  F(e^{it})\big)}
{(e^{-it}-\overline{z})^{\frac{\alpha}{2}}(1-ze^{-it})^{\frac{\alpha}{2}+1}} \,dt  \\
=&
-\frac{ \mathcal{C}_{\alpha}(1-| z |^2)^\alpha}{  \pi}
\int_{0}^{2\pi}\dot F\left(e^{i(t+\theta)}\right)
 \frac{r\sin t }{|1-re^{i t }|^{\alpha+2}}\,dt\\
&-\frac{\alpha \mathcal{C}_{\alpha}(1-| z |^2)^\alpha}{ 2\pi}
\int_{0}^{2\pi} F\left(e^{i(t+\theta)}\right)
 \frac{1-r\cos t }{|1-re^{i t }|^{\alpha+2}} \,dt,
\end{split}
\end{align*}
and so, we get
\begin{align}\label{eq-3.18}
\begin{split}
|\mathscr{J}_{\alpha,2} (re^{i\theta}) |
\leq  &\frac{\mathcal{C}_{\alpha} (1-r^2)^\alpha}{ \pi} \mathscr{ A} _{\alpha,1}(r,\theta)
 +
\frac{|\alpha |\mathcal{C}_{\alpha}(1-r^2)^\alpha}{2\pi}\mathscr{ A} _{\alpha,2}(r,\theta),
\end{split}
\end{align}
where
\begin{equation}\label{eq-3.19}
\mathscr{ A} _{\alpha,1}(r,\theta)
=\int_{0}^{2\pi} \left|\dot F\left(e^{i(t+\theta)}\right)\right|\frac{ r|\sin t|}{|1-re^{it}|^{\alpha+2}}\,dt
\end{equation}
and
\begin{equation}\label{eq-3.20}
\mathscr{ A}_{\alpha,2}(r,\theta)
=\int_{0}^{2\pi} \left|F\left(e^{i(t+\theta)}\right)\right| \frac{|1-r\cos t|}{| 1-re^{i t}|^{\alpha+2}}\,dt.
\end{equation}
It follows from  $\dot F\in \mathcal{L}^{\infty}(\mathbb{T})$  and
 $$\max\{|r\sin t|,|1-r\cos t|\}\leq |1-re^{it}|$$
 that
\begin{align*}
\begin{split}
\;\;\;\;\;|\mathscr{J}_{\alpha,2}(z)|
&\leq    \frac{\mathcal{C}_{\alpha}  }{2\pi}\big(2 \| \dot{F}\|_{\mathcal{L}^{\infty}(\mathbb{T})}  +|\alpha|\cdot
 \| F\|_{\mathcal{L}^{\infty}(\mathbb{T})}  \big )
\int_{0}^{2\pi}\frac{(1-r^{2})^{\alpha}}{ |1-re^{it}|^{\alpha+1}}dt.
\end{split}
\end{align*}
Then for $\alpha\in(0,\infty)$, Lemma I guarantees that
\begin{align*}
\begin{split}
 |\mathscr{J}_{\alpha,2}(z)|
&\leq   \mathcal{C}_{\alpha}\big(2 \| \dot{F}\|_{\mathcal{L}^{\infty}(\mathbb{T})}  +\alpha \| F\|_{\mathcal{L}^{\infty}(\mathbb{T})}   \big)
\frac{\Gamma(\alpha)}{\Gamma^{2}(\frac{\alpha+1}{2})},
\end{split}
\end{align*}
which shows that
\begin{align}\label{eq-3.21}
\begin{split}
\|\mathscr{J}_{\alpha,2}\|_{\mathcal{H}_{\mathcal{G}}^{\infty}(\mathbb{D})}
\lesssim
\|F\|_{\mathcal{L}^{\infty}(\mathbb{T})}+\|\dot{F}\|_{\mathcal{L}^{\infty}(\mathbb{T})}.
\end{split}
\end{align}

\bca\label{case-4}  Suppose that $p\in[1,\infty)$.
\eca

 Let $z=re^{i\theta}\in \mathbb{D}$.
 For  $\alpha\in (-1,\infty)$,  we infer from (\ref{eq-3.15}) that
\begin{equation}\label{eq-3.22}
| \mathscr{J}_{\alpha,1}(re^{i\theta})|\leq\frac{|\alpha| \mathcal{C}_{\alpha}}{2\pi}\int_{0}^{2\pi}
\frac{(1-r^2)^{\alpha+1}}{| 1-re^{i t}|^{\alpha+2}}\left| F\left(e^{i(t+\theta)}\right)\right|\,dt.
\end{equation}
Then \eqref{eq-3.22}  yields that for any $p \in[1,\infty)$,
$$
|\mathscr{J}_{\alpha,1}(re^{i\theta})| ^p
\leq  \left(\frac{|\alpha| \mathcal{C}_{\alpha}}{2\pi}\right)^p I_{\alpha,1}^p(r)
\left(\int_{0}^{2\pi}\frac{(1-r^2)^{\alpha+1}\left| F\left(e^{i(t+\theta)}\right)\right|}
{| 1-re^{i t}|^{\alpha+2}I_{\alpha,1}(r)}
\,dt\right)^p,
 $$
 where
$I_{\alpha,1}$ is defined in  \eqref{eq-3.2}.
Further, by applying Jensen's inequality,
we have
\begin{equation*}
| \mathscr{J}_{\alpha,1}(re^{i\theta})| ^p
\leq \left(\frac{|\alpha| \mathcal{C}_{\alpha}}{2\pi}\right)^p I_{\alpha,1}^{p-1}(r)
\int_{0}^{2\pi}\frac{(1-r^2)^{\alpha+1}}{| 1-re^{i t}|^{\alpha+2}}\left| F\left(e^{i(t+\theta)}\right)\right|^p\,dt.
\end{equation*}
Therefore,
\begin{align*}
\begin{split}
&H_{p}^{p}\big(\mathscr{J}_{\alpha,1}\big)
= \sup_{r\in[0,1)}\frac{1}{2\pi}\int_{0}^{2\pi} | \mathscr{J}_{\alpha,1}(re^{i\theta})|^p\,d\theta\\
&\leq\frac{\left(|\alpha| \mathcal{C}_{\alpha}\right)^p}{(2\pi)^{p+1}} \sup_{r\in[0,1)}I_{\alpha,1}^{p-1}(r)\int_{0}^{2\pi}
\int_{0}^{2\pi}\frac{(1-r^2)^{\alpha+1}}{|  1-re^{i t }|^{\alpha+2}}\left| F\left(e^{i(t+\theta)}\right)\right|^p\,dt\,d\theta\\
&\leq \left(\frac{|\alpha| \mathcal{C}_{\alpha}}{2\pi}\right)^p\| F \|_{\mathcal{L}^{\infty}(\mathbb{T})}^{p}\sup_{r\in[0,1)}I_{\alpha,1}^{p}(r),
\end{split}
\end{align*}
where the last inequality follows from \eqref{eq-3.2} and the assumption that $F$ is continuous on $\mathbb{T}$.
Thus for $\alpha\in(-1,\infty)$ and $p\in [1,\infty)$, we infer from \eqref{eq-1.5} and \eqref{eq-3.3} that
\begin{align}\label{eq-3.23}
\begin{split}
\|\mathscr{J}_{\alpha,1}\|^{p}_{\mathcal{H}_{\mathcal{G}}^{p}(\mathbb{D})}
&\leq |\alpha|^p\| F \|_{\mathcal{L}^{\infty}(\mathbb{T})}^{p}.
\end{split}
\end{align}
This, together with \eqref{eq-3.17}, implies that for $\alpha\in(-1,\infty)$ and $p\in [1,\infty)$
or for $\alpha\in(0,\infty)$ and $p=\infty$,
\begin{align*}
\begin{split}
\|\mathscr{J}_{\alpha,1}\|_{\mathcal{H}_{\mathcal{G}}^{p}(\mathbb{D})}\lesssim\|F\|_{\mathcal{L}^{\infty}(\mathbb{T})}.
\end{split}
\end{align*}
 Hence the first statement in the lemma is true.

For the proofs of the second and the third statements, we need some preparation.
For $\alpha\in(-1,\infty)$ and $p\in [1,\infty)$, we infer from Lemma H and \eqref{eq-3.18} that
\begin{align}\label{eq-3.24}
\begin{split}
|\mathscr{J}_{\alpha,2} (re^{i\theta}) |^p
\leq&\frac{ (\mathcal{C}_{\alpha})^p(1-r^2)^{\alpha p}}{2 \pi^p}
\left(2^{p}\mathscr{A}_{\alpha,1}^p(r,\theta)+|\alpha| ^{p}\mathscr{A}_{\alpha,2}^p(r,\theta) \right).
%
\end{split}
\end{align}

In order to estimate
$\mathscr{A}_{\alpha,1}^p(r,\theta) $ and $\mathscr{A}_{\alpha,2}^p(r,\theta) $, let
\begin{equation}\label{eq-3.25}
I_{\alpha,2}(r):= \int_{0}^{2\pi}\frac{ r |\sin t|}{| 1-re^{it}|^{\alpha+2}}\,dt
\;\;\text{and}\;\;
I_{\alpha,3}(r):= \int_{0}^{2\pi}\frac{|1-r\cos t|}{| 1-re^{it}|^{\alpha+2}}\,dt.
\end{equation}
Then it follows from \eqref{eq-3.19} and Jensen's inequality that
\begin{align}\label{eq-3.26}
\begin{split}
\mathscr{A}_{\alpha,1}^p(r,\theta) &=I_{\alpha,2}^p(r)\left(\int_{0}^{2\pi}\left|\dot F\left(e^{i(t+\theta)}\right)\right| \frac{ r|\sin t|}{| 1-re^{it}|^{\alpha+2}I_{\alpha,2}(r)}
\,dt\right)^p\\
&\leq I_{\alpha,2}^{p-1}(r)\int_{0}^{2\pi}\left|\dot F\left(e^{i(t+\theta)}\right)\right|^{p} \frac{ r|\sin t|}{| 1-re^{it}|^{\alpha+2}}\,dt.
\end{split}
\end{align}
Similarly, we get from \eqref{eq-3.20} that
\begin{align}\label{eq-3.27}
\begin{split}
 \mathscr{A}_{\alpha,2}^p(r,\theta)
\leq I_{\alpha,3}^{p-1}(r)\int_{0}^{2\pi}\left|  F\left(e^{i(t+\theta)}\right)\right|^p\frac{ |1-r\cos t|}{| 1-re^{it }|^{\alpha+2}}\,dt.
\end{split}
\end{align}

Substituting \eqref{eq-3.26} and (\ref{eq-3.27}) into \eqref{eq-3.24} gives
\begin{align*}
\begin{split}
|\mathscr{J}_{\alpha,2} (re^{i\theta}) |^p
\leq&\frac{2^{p-1}(\mathcal{C}_{\alpha})^p(1-r^2)^{\alpha p}}{ \pi^p}
 I_{\alpha,2}^{p-1}(r)\int_{0}^{2\pi }\left|  \dot F\left(e^{i(t+\theta)}\right)\right|^p\frac{r|\sin t|}{| 1-re^{it}|^{\alpha+2}}\,dt \\
& +  \frac{ |\alpha|^{p} (\mathcal{C}_{\alpha})^p(1-r^2)^{\alpha p}}{2 \pi^p}I_{\alpha,3}^{p-1}(r) \int_{0}^{2\pi}\left|  F\left(e^{i(t+\theta)}\right)\right|^p\frac{ |1-r\cos t|}{| 1-re^{it}|^{\alpha+2}}\,dt .
\end{split}
\end{align*}

Since \eqref{eq-1.7} and the assumption of $\dot F \in \mathcal{L}^{p}(\mathbb{T})$ ensure that
$$
\int_{0}^{2\pi}\int_{0}^{2\pi}\left| \dot F\left(e^{i(t+\theta)}\right)\right|^p \frac{r|\sin t|}{| 1-re^{it}|^{\alpha+2}}\, dt d\theta
\leq\frac{4\pi^{2}}{(1-r)^{\alpha+2}}\|  \dot F\|^{p}_{\mathcal{L}^{p}(\mathbb{T})}
<\infty
$$
and
$$
\int_{0}^{2\pi}\int_{0}^{2\pi}\left|  F\left(e^{i(t+\theta)}\right)\right|^p\frac{ |1-r\cos t|}{| 1-re^{it}|^{\alpha+2}}\,dt d\theta
\leq\frac{8\pi^{2}}{(1-r)^{\alpha+2}} \|F \|^{p}_{\mathcal{L}^{\infty}(\mathbb{T})}
<\infty,
$$
we know from Fubini's Theorem that
\begin{align*}
\begin{split}
\int_{0}^{2\pi}&|\mathscr{J}_{\alpha,2}  (re^{i\theta}) |^p\,d\theta\\
\leq&\frac{2^{p-1}(\mathcal{C}_{\alpha})^p(1-r^2)^{\alpha p}}{ \pi^p}
  I_{\alpha,2}^{p-1}(r)\int_{0}^{2\pi}\frac{r|\sin t|}{|1-re^{it}|^{\alpha+2}}\,dt
\int_{0}^{2\pi}\left| \dot F\left(e^{i(t+\theta)}\right)\right|^p\,d\theta \\
& +   \frac{|\alpha|^{p}(\mathcal{C}_{\alpha})^p(1-r^2)^{\alpha p}}{2 \pi^p}  I_{\alpha,3}^{p-1}(r) \int_{0}^{2\pi}\frac{ |1-r\cos t|}{| 1-re^{it}|^{\alpha+2}}\,dt
\int_{0}^{2\pi}\left|  F\left(e^{i(t+\theta)}\right)\right|^p\,d\theta  .
\end{split}
\end{align*}
This, together with \eqref{eq-1.7} and \eqref{eq-3.25}, implies that
\begin{align*}
\begin{split}
 \int_{0}^{2\pi}|\mathscr{J}_{\alpha,2}  (re^{i\theta}) |^p\,d\theta
\leq\frac{(\mathcal{C}_{\alpha})^p(1-r^2)^{\alpha p} }{\pi^{p-1}}\left( 2^{p}I_{\alpha,2}^{p }(r)\|\dot F \|_{\mathcal{L}^p(\mathbb{T})}^{p}
+ |\alpha|^{p}I_{\alpha,3}^{p }(r) \| F \|_{\mathcal{L}^{\infty}(\mathbb{T})}^{p}\right) .
\end{split}
\end{align*}
Then for $\alpha\in(-1,\infty)$ and $p\in[1,\infty)$,
\begin{align}\label{eq-3.28}
\begin{split}
 H_{p}^{p}\big(\mathscr{J}_{\alpha,2}\big)
 \leq
     \frac{(\mathcal{C}_{\alpha})^p }{2\pi^{p}}   \sup_{r\in[0,1)}
     \left( 2^{p}\|\dot F \|_{\mathcal{L}^p(\mathbb{T})}^{p} \mathscr{B}_{\alpha,1}^{p}(r)
+ |\alpha|^{p}\| F \|_{\mathcal{L}^{\infty}(\mathbb{T})}^{p}\mathscr{B}_{\alpha,2}^{p}(r)\right)
\end{split}
\end{align}
 and
\begin{align}\label{eq-3.29}
\begin{split}
 B_{p}^{p}\big(\mathscr{J}_{\alpha,2}\big)
 \leq\frac{(  \mathcal{C}_{\alpha})^{p}}{\pi^{p}}
 \left(2^{p} \|\dot F\|^{p}_{\mathcal{L}^{p} (\mathbb{T})} \mathscr{C}_{\alpha,p}
 +|\alpha  |^{p}   \|F\|^{p}_{\mathcal{L}^{\infty} (\mathbb{T})} \mathscr{D}_{\alpha,p}
\right) ,
\end{split}
\end{align}
where
\be\label{eq-3.30}
\mathscr{B}_{\alpha,1}(r)
 := (1-r^2)^{\alpha}  I_{\alpha,2} (r),\quad
 \mathscr{B}_{\alpha,2}(r)
 := (1-r^2)^{\alpha}   I_{\alpha,3} (r),
\ee
\be\label{eq-3.31}
\mathscr{C}_{\alpha,p}
 := \int_{0}^{1}r(1-r^{2})^{\alpha p}I^{p}_{\alpha,2}(r)\,dr
 \; \;\text{and}\;\;
 \mathscr{D}_{\alpha,p}
 := \int_{0}^{1}r(1-r^{2})^{\alpha p}I^{p}_{\alpha,3}(r)\,dr.
\ee

Still, we need the following estimates of the quantities $\mathscr{B}_{\alpha,1}^{p}$, $\mathscr{B}_{\alpha,2}^{p}$,
$\mathscr{C}_{\alpha,p}$ and $\mathscr{D}_{\alpha,p}$.
\medskip

\noindent
$(\mathfrak{A})$ {\bf Estimate on  $\mathscr{B}_{\alpha,1}^{p}$}. It follows from \eqref{eq-3.25} and \eqref{eq-3.30}  that
\begin{align}\label{eq-3.32}
\begin{split}
\mathscr{B}_{\alpha,1} (r)
&=\left\{
    \begin{array}{ll}
      \frac{2}{\alpha} \big( (1+r)^{\alpha}- (1-r)^{\alpha}\big), & \hbox{if\;\;$\alpha\in(-1,0)\cup(0,\infty)$,} \\
      2\log\frac{1+r}{1-r}, & \hbox{if\;\;$\alpha=0$.}
    \end{array}
  \right.
\end{split}
\end{align}
Then for  $\alpha\in(0,\infty)$ and $p\in[1,\infty)$,
\be\label{eq-3.33}
0\leq \mathscr{B}_{\alpha,1}^{p} (r)\leq \frac{2^{p}}{\alpha^{p}} (1+r)^{\alpha p} < \frac{2^{(1+\alpha)p}}{  \alpha^{ p}}.
\ee

\noindent
$(\mathfrak{B})$ {\bf Estimate on  $\mathscr{C}_{\alpha,p}$}.
 For the case when $\alpha=0$ and $p\in[1,\infty)$, we infer from \eqref{eq-3.30}$-$\eqref{eq-3.32} that
\begin{align*}
\begin{split}
\mathscr{C}_{\alpha,p}
=\int_{0}^{1} r \mathscr{B}_{\alpha,1}^{p} (r)dr
\leq2^{p}\int_{0}^{1}  \log^{p}\frac{2}{1-r}dr
= 2^{p} \left(\int_{0}^{ \frac{1}{2}}\log^{p}\frac{2}{s}\,ds +\int_{ \frac{1}{2}}^{1}\log^{p}\frac{2}{s}\,ds\right),
\end{split}
\end{align*}
where $s=1-r\in(0,1]$.
Since
$$
\int_{ \frac{1}{2}}^{1}\log^{p}\frac{2}{s}\,ds\leq  \frac{ \log^{p}4}{2},
$$
and since for $s\in(0,\frac{1}{2}]$, the fact $\log \frac{2}{s}\leq 2\log\frac{1}{s}$ leads to
$$
\int_{0}^{ \frac{1}{2}}\log^{p}\frac{2}{s}\,ds
\leq2^{p}\int_{0}^{ \frac{1}{2}}\log^{p}\frac{1}{s}\,ds
\leq2^{p}\int_{0}^{1}\log^{p}\frac{1}{s}\,ds
=2^{p}\Gamma(1+p),
$$
we have that for $\alpha=0$ and $p\in[1,\infty)$,
\begin{align}\label{eq-3.34}
\begin{split}
0\leq
\mathscr{C}_{\alpha,p}
 \leq  2^{p-1}  \log^{p}4+4^{p}\Gamma(1+p).
\end{split}
\end{align}

For the case when $\alpha\in(-1,0)$ and
$p\in\big[1,-\frac{1}{\alpha}\big)$, since $$|r\sin t|\leq |1-re^{it}|,$$
we infer from \eqref{eq-3.25}  and \eqref{eq-3.31} that
\begin{equation}\label{eq-3.35}
\begin{split}
\mathscr{C}_{\alpha,p}
\leq&  \int_{0}^{1} r\left(  \int_{0}^{2\pi} \frac{(1-r^{2})^{\alpha  }}{|1-re^{it}|^{\alpha+1}}dt\right)^{p}dr \\
=&\int_{0}^{\frac{1}{2}} r \left(  \int_{0}^{2\pi} \frac{(1-r^{2})^{\alpha  }}{|1-re^{it}|^{\alpha+1}}dt\right)^{p}dr
  +  \int_{\frac{1}{2}}^{1} r \left(  \int_{0}^{2\pi} \frac{(1-r^{2})^{\alpha  }}{|1-re^{it}|^{\alpha+1}}dt\right)^{p}dr.
\end{split}
\end{equation}

Moreover, we know from the fact ``$\frac{ r(1+r)^{\alpha p}}{(1-r)^{p}}\leq 2^{p-1}$ for $r\in[0,\frac{1}{2}]$" that
\begin{equation}\label{eq-3.36}
\begin{split}
 \int_{0}^{\frac{1}{2}} r\left(  \int_{0}^{2\pi} \frac{(1-r^{2})^{\alpha  }}{|1-re^{it}|^{\alpha+1}}dt\right)^{p}dr
\leq &\int_{0}^{\frac{1}{2}} r \left(  \int_{0}^{2\pi} \frac{(1-r^{2})^{\alpha  }}{(1-r)^{\alpha+1}}dt\right)^{p}dr\\
\leq& 2^{2p-2}\pi^{p}.
\end{split}
\end{equation}
 Since
\begin{equation*}
\begin{split}
\int_{\frac{1}{2}}^{1} r(1-r^{2})^{\alpha p} dr
\leq \frac{3^{\alpha p}}{2^{\alpha p}}\int_{\frac{1}{2}}^{1}(1-r )^{\alpha p} dr
\leq\frac{3^{\alpha p}}{2^{\alpha p}}\int_{0}^{1}s^{\alpha p} ds
=\frac{3^{\alpha p}}{2^{\alpha p} (\alpha p+1)},
\end{split}
\end{equation*}
 by Lemma \ref{lem-3.2}, we have
\begin{equation}\label{eq-3.37}
  \begin{split}
  \int_{ \frac{1}{2}}^{1} r\left(  \int_{0}^{2\pi} \frac{(1-r^{2})^{\alpha  }}{|1-re^{it}|^{\alpha+1}}dt\right)^{p}dr
  \leq \frac{1}{\alpha p+1}
\left(  \frac{3^{\frac{3\alpha+1}{2}} }{2^{ 2\alpha-1 } }  \frac{\Gamma(-\alpha) \Gamma(\frac{1}{2})}{\Gamma(\frac{1}{2}-\alpha)} \right)^{p}.
\end{split}
\end{equation}

Substituting  \eqref{eq-3.36} and   \eqref{eq-3.37} into \eqref{eq-3.35}  yields that for $\alpha\in(-1,0)$ and
$p\in\big[1,-\frac{1}{\alpha}\big)$,
\begin{equation}\label{eq-3.38}
\begin{split}
\mathscr{C}_{\alpha,p}
\leq  \int_{0}^{1} r\left(  \int_{0}^{2\pi} \frac{(1-r^{2})^{\alpha  }}{|1-re^{it}|^{\alpha+1}}dt\right)^{p}dr
\leq  M_{3},
\end{split}
\end{equation}
where
\begin{equation*}
\begin{split}
M_{3} = 2^{2p-2}\pi^{p}+\frac{1}{\alpha p+1}
\left(  \frac{3^{\frac{3\alpha+1}{2}} }{2^{ 2\alpha-1 } }  \frac{\Gamma(-\alpha) \Gamma(\frac{1}{2})}{\Gamma(\frac{1}{2}-\alpha)} \right)^{p}.
\end{split}
\end{equation*}

\noindent $(\mathfrak{C})$ {\bf Estimate on  $\mathscr{B}_{\alpha,2}^{p}$}.
Since  $$ |1-r\cos t| \leq |1-re^{i t}|,$$   by \eqref{eq-3.25} and \eqref{eq-3.30}, we find that
for $\alpha\in(-1,\infty)$,
\begin{equation}\label{eq-3.39}
\mathscr{B}_{\alpha,2} (r)
\leq \int_{0}^{2\pi}\frac{(1-r^{2})^{\alpha}}{|1-re^{i t}|^{\alpha+1}}\,dt.
\end{equation}
Then we deduce from Lemma I and the fact  ``$a\Gamma(a)=\Gamma(a+1)$ for $a>0$" that
\begin{equation*}
\alpha\mathscr{B}_{\alpha,2} (r)
\leq\left\{
    \begin{array}{ll}
     \frac{2\pi\Gamma(\alpha+1 )}{\Gamma^2(\frac{\alpha+1}{2})}, & \hbox{if\;\;$\alpha\in(0,\infty)$;} \\
      0, & \hbox{if\;\;$\alpha=0$.}
    \end{array}
  \right.
\end{equation*}
Thus for $\alpha\in[0,\infty)$ and $p\in[1,\infty)$, we have
\begin{equation}\label{eq-3.40}
0\leq \alpha^{p} \mathscr{B}_{\alpha,2}^{p}(r)\leq \frac{2^{p}\pi^{p}\Gamma^{p}(\alpha+1)}{\Gamma^{2p}(\frac{\alpha+1}{2})}.
\end{equation}

\noindent $(\mathfrak{D})$ {\bf Estimate on $\mathscr{D}_{\alpha,p}$}.    If $\alpha\in(-1,0)$ and
$p\in\big[1,-\frac{1}{\alpha}\big)$,  by
\eqref{eq-3.30}, \eqref{eq-3.31} and \eqref{eq-3.39}, we arrive at the following inequality
\begin{equation*}
 \mathscr{D}_{\alpha,p}
 = \int_{0}^{1} r \mathscr{B}^{p}_{\alpha,2}(r)dr
\leq \int_{0}^{1} r\left( \int_{0}^{2\pi} \frac{(1-r^{2})^{\alpha  }}{|1-re^{it}|^{\alpha+1}}dt\right)^{p}dr.
\end{equation*}
Then it follows from \eqref{eq-3.38} that
\begin{equation}\label{eq-3.41}
 \mathscr{D}_{\alpha,p}
\leq \int_{0}^{1} r\left( \int_{0}^{2\pi} \frac{(1-r^{2})^{\alpha  }}{|1-re^{it}|^{\alpha+1}}dt\right)^{p}dr
\leq M_{3}.
\end{equation}

Now, we are ready to prove the rest two statements of the lemma. Firstly, we conclude from  \eqref{eq-3.28}, \eqref{eq-3.33} and \eqref{eq-3.40} that for $\alpha\in(0,\infty)$ and $p\in[1,\infty)$,
\begin{align*}
\begin{split}
 H_{p}^{p}\big(\mathscr{J}_{\alpha,2}\big)
& \leq
     \frac{(\mathcal{C}_{\alpha})^{p}}{2\pi^{p }}
\left(\frac{2^{(2+\alpha)p}}{\alpha^{p}}\|\dot F \|_{\mathcal{L}^p(\mathbb{T})}^{p}
+ \frac{2^{p}\pi^p \Gamma^{p}(\alpha+1)}{\Gamma^{2p}(\frac{\alpha+1}{2})} \| F \|_{\mathcal{L}^{\infty}(\mathbb{T})}^{p}\right),
\end{split}
\end{align*}
which, together with \eqref{eq-3.21}, indicates that for $\alpha\in(0,\infty)$ and $p\in[1,\infty]$,
$$\|\mathscr{J}_{\alpha,2}\|_{\mathcal{H}_{\mathcal{G}}^{p}(\mathbb{D})}
\lesssim\|F\|_{\mathcal{L}^{\infty}(\mathbb{T})}+\|\dot{F}\|_{\mathcal{L}^{p}(\mathbb{T})}.$$
Hence the proof of the second statement of the lemma is complete.

Secondly, it follows from \eqref{eq-3.29} and \eqref{eq-3.34}  that for $\alpha=0$ and $p\in[1,\infty)$,
\begin{align}\label{eq-3.42}
\begin{split}
 B_{p}^{p}\big(\mathscr{J}_{\alpha,2}\big)
&\leq \frac{ (\mathcal{C}_{\alpha})^p }{\pi^{p}}
\big(2^{2p-1}\log^{p}4+8^{p}\Gamma(1+p)   \big)
\|\dot F \|_{\mathcal{L}^p(\mathbb{T})}^{p},
\end{split}
\end{align}
and from \eqref{eq-3.29}, \eqref{eq-3.38} and \eqref{eq-3.41} that for $\alpha\in(-1,0)$ and
$p\in\big[1,-\frac{1}{\alpha}\big)$,
\begin{align}\label{eq-3.43}
\begin{split}
 B_{p}^{p}\big(\mathscr{J}_{\alpha,2}\big)
\leq&\frac{  (\mathcal{C}_{\alpha})^{p}M_{3} }{\pi^{p}}
\left(2^{p} \|\dot F\|^{p}_{\mathcal{L}^{p} (\mathbb{T}) }
+|\alpha|^{p}\|F\|^{p}_{\mathcal{L}^{\infty} (\mathbb{T})}\right)<\infty.
\end{split}
\end{align}
These show that if $\alpha=0$ and $p\in[1,\infty)$ or if $\alpha\in(-1,0)$ and
$p\in\big[1,-\frac{1}{\alpha}\big)$, then $$\|\mathscr{J}_{\alpha,2}\|_{\mathcal{B}_{\mathcal{G}}^{p}(\mathbb{D})}
\lesssim\|F\|_{\mathcal{L}^{\infty}(\mathbb{T})}+\|\dot{F}\|_{\mathcal{L}^{p}(\mathbb{T})}.$$
This demonstrates that the third statement of the lemma is true.
\epf

Based on Lemma \ref{lem-3.3}, we are going to prove the following lemma.
\begin{lem}\label{lem-3.4}
Under the assumptions of Lemma \ref{lem-3.1}, the   following two statements are true:
\begin{enumerate}
  \item For all $\alpha\in(0,\infty)$ and $p\in[1,\infty]$,
$$
  \|\partial_{r}f \|_{\mathcal{H}_{\mathcal{G}}^{p}(\mathbb{D})}\lesssim \|F\|_{\mathcal{L}^{\infty}(\mathbb{T})}+\|\dot{F}\|_{\mathcal{L}^{p}(\mathbb{T})};
  $$
  \item For $\alpha=0$ and  $p=1$ or  for all $\alpha\in(-1,0)$ and
$p\in\big[1,-\frac{1}{\alpha}\big)$,
$$
\|\partial_{r}f \|_{\mathcal{B}_{\mathcal{G}}^{p}(\mathbb{D})}\lesssim \|F\|_{\mathcal{L}^{\infty}(\mathbb{T})}+\|\dot{F}\|_{\mathcal{L}^{p}(\mathbb{T})}.
$$
\end{enumerate}

\end{lem}
 \bpf For $z=re^{i\theta}\in \mathbb{D}$ and $t\in[0,2\pi]$, it follows from  \eqref{eq-1.4}   that
 $$
\frac{\partial  \mathcal{K}_{\alpha}(ze^{-it})}{\partial z}F(e^{it})
= \mathcal{C}_{\alpha}(1-| z |^2)^{\alpha} e^{-it} \frac{(\alpha+2)(1-\overline{z}e^{it})+ \alpha \overline{z}(z-e^{it})}
{2(1-ze^{-it})^{\frac{\alpha}{2}+2}(1-\overline{z}e^{it})^{\frac{\alpha+2}{2}}}F(e^{it})
 $$
 and
 $$
 \frac{\partial  \mathcal{K}_{\alpha}(ze^{-it})}{\partial \overline{z}}F(e^{it})
  = \mathcal{C}_{\alpha}(1-| z |^2)^{\alpha} e^{it}
\frac{(\alpha+2)(1-ze^{-it})+ \alpha z(\overline{z}-e^{-it})}
{2(1-ze^{-it})^{\frac{\alpha+2}{2}}(1-\overline{z}e^{it})^{\frac{\alpha}{2}+2}}F(e^{it}),
$$
where $\mathcal{C}_{\alpha}$ is the constant from (\ref{eq-1.5}).
Since the mappings
$$(z,t)\mapsto \mathcal{K}_{\alpha}(ze^{-it})F(e^{it}), \;
(z,t)\mapsto  \frac{\partial  \mathcal{K}_{\alpha}(ze^{-it})}{\partial z}F(e^{it})
\;\text{and}\;
(z,t)\mapsto  \frac{\partial  \mathcal{K}_{\alpha}(ze^{-it})}{\partial \overline{z}}F(e^{it})
 $$
are continuous in $\overline{\mathbb{D}}_{\rho_{0}}\times[0,2\pi]$ for any $\rho_{0}\in(0,1)$,
 we obtain from \eqref{eq-1.3} that
 \begin{align*}
\begin{split}
\partial_{z}f(z)
=&\frac{1}{2\pi}\int_{0}^{2\pi}\frac{\partial  \mathcal{K}_{\alpha}(ze^{-it})}{\partial z}F(e^{it})\,dt\\
=&\frac{\mathcal{C}_{\alpha}(1-| z |^2)^{\alpha}}{2\pi}\int_{0}^{2\pi}\left(\frac{(\alpha+2)e^{it}}
{2(e^{-it}-\overline{z})^{\frac{\alpha}{2}}(e^{it}-z)^{\frac{\alpha}{2}+2}}-\frac{ \alpha \overline{z}}{2| e^{it}-z | ^{\alpha+2} }\right)F(e^{it})\,dt,
\end{split}
\end{align*}

\begin{align*}
 \begin{split}
\partial_{\overline{z}}f(z)
=&\frac{1}{2\pi}\int_{0}^{2\pi}
  \frac{\partial  \mathcal{K}_{\alpha}(ze^{-it})}{\partial \overline{z}}F(e^{it}) \,dt\\
=&\frac{\mathcal{C}_{\alpha}(1-| z |^2)^{\alpha}}{2\pi}\int_{0}^{2\pi}\left(\frac{(\alpha+2)e^{-it}}
{2(e^{it}-z)^{\frac{\alpha}{2}}(e^{-it}-\overline{z})^{\frac{\alpha}{2}+2}}-\frac{ \alpha z}{2| e^{it}-z | ^{\alpha+2} }\right)F(e^{it})\,dt
 \end{split}
\end{align*}
and $\partial_{z}f,\partial_{\overline{z}}f\in \mathcal{C}(\mathbb{D})$.
Therefore,
\begin{align*}
z\partial_{z}f(z)+\overline{z}\partial_{\overline{z}}f(z)
=&\frac{\mathcal{C}_{\alpha}(1-| z |^2)^{\alpha}}{2\pi}\int_{0}^{2\pi}\left(\frac{(\alpha+2)ze^{it}}
{2(e^{-it}-\overline{z})^{\frac{\alpha}{2}}(e^{it}-z)^{\frac{\alpha}{2}+2}}\right.\\
&\left.+\frac{(\alpha+2)\overline{z}e^{-it}}
{2(e^{it}-z)^{\frac{\alpha}{2}}(e^{-it}-\overline{z})^{\frac{\alpha}{2}+2}}
-\frac{\alpha | z| ^2}{| e^{it}-z | ^{\alpha+2} }\right)F(e^{it})\,dt.
\end{align*}

Moreover, elementary calculations give that
\begin{align*}\label{eq-3.2}
 \begin{split}
 \frac{(\alpha+2)ze^{it}}
{2(e^{-it}-\overline{z})^{\frac{\alpha}{2}}(e^{it}-z)^{\frac{\alpha}{2}+2}}
=\frac{\alpha}{2| e^{it}-z|^{\alpha+2}}
  +  e^{-  \frac{i\alpha t}{2}}
\left(\frac{\partial}{\partial{t}}\frac{i}{(e^{-it}-\overline{z})^{\frac{\alpha}{2}}(1-ze^{-it})^{\frac{\alpha}{2}+1}}\right)
 \end{split}
\end{align*}
and
\begin{align*}
 \begin{split}
 \frac{(\alpha+2)\overline{z}e^{-it}}
{2(e^{it}-z)^{\frac{\alpha}{2}}(e^{-it}-\overline{z})^{\frac{\alpha}{2}+2}}
=\frac{\alpha}{2| e^{it}-z|^{\alpha+2}}  -  e^{  \frac{i\alpha t}{2}} \left(\frac{\partial}{\partial{t}}\frac{i}{(e^{it}-z)^{\frac{\alpha}{2}}(1-\overline{z}e^{it})^{\frac{\alpha}{2}+1}}\right).
 \end{split}
\end{align*}

Combining the above three equalities yields that
\begin{equation}\label{eq-3.44}
  r\partial_{r}f(z)=z\partial_{z}f(z)+\overline{z}\partial_{\overline{z}}f(z)=\mathscr{J}_{\alpha,1}(z)+\mathscr{J}_{\alpha,2}(z),
\end{equation}
where $\mathscr{J}_{\alpha,1}$ and $\mathscr{J}_{\alpha,2}$ are defined in \eqref{eq-3.15} and \eqref{eq-3.16}, respectively.

In the following, we divide the discussions into two cases.
\bca\label{case-5}  Suppose that $p\in[1,\infty)$.
\eca

Since \eqref{eq-3.44} and the continuity of $\partial_{z}f$, $\partial_{\overline{z}}f$ imply that $\partial_{r}f$ is continuous in
$ \mathbb{D} $,
we see that
\begin{align}\label{eq-3.45}
\begin{split}
\sup_{r\in[0,\frac{1}{2}]} \frac{1}{2\pi }\int_{0}^{2\pi}|\partial_{r}f(re^{i\theta})|^{p}d\theta
\leq\ M_{4}^{p}
<\infty
\end{split}
\end{align}
and
\begin{align}\label{eq-3.46}
\begin{split}
 \int_{ \overline{\mathbb{D}}_{\frac{1}{2}}} |\partial_{r}f(z)|^{p}d\sigma(z)
=\frac{1}{\pi}\int_{0}^{\frac{1}{2}}\int_{0}^{2\pi}\big|\partial_{r}f(re^{i\theta})\big|^{p}rd\theta dr
\leq M_{4}^{p},
\end{split}
\end{align}
where
\be\label{eq-3.47}
M_{4} =\max\left\{\| D_{f}(z)\| :z\in \overline{\mathbb{D}}_{\frac{1}{2}}\right\}.
\ee
For $z=re^{i\theta}\in \mathbb{D}\backslash \overline{\mathbb{D}}_{\frac{1}{2}}$, $\alpha\in(-1,\infty)$ and  $p\in [1,\infty)$,
by Lemma H and (\ref{eq-3.44}), we have that
\begin{align}\label{eq-3.48}
\begin{split}
| \partial_{r}f(z)| ^p
&\leq  \frac{1}{r^{p}}\left( | \mathscr{J}_{\alpha,1}(z)| + | \mathscr{J}_{\alpha,2}(z)|  \right)^p
\leq  2^{p}\left( | \mathscr{J}_{\alpha,1}(z)| + | \mathscr{J}_{\alpha,2}(z)|   \right)^p \\
&\leq 2^{2p-1}
\left( | \mathscr{J}_{\alpha,1}(z)|^p + | \mathscr{J}_{\alpha,2}(z)|^p \right).
\end{split}
\end{align}

Now, we are ready to estimate the quantities $\|\partial_{r}f\|_{\mathcal{H}_{\mathcal{G}}^{p}(\mathbb{D})} $ and $\|\partial_{r}f\|_{\mathcal{B}_{\mathcal{G}}^{p}(\mathbb{D})} $.
\medskip

\noindent
 (I)  {\bf Estimate on  $\|\partial_{r}f\|_{\mathcal{H}_{\mathcal{G}}^{p}(\mathbb{D})} $}.
For $\alpha\in(0,\infty)$ and $p\in[1,\infty)$, we have   that
\begin{align*}
\begin{split}
H_{p}^{p}\big(\partial_{r}f \big)
\leq& M_{4}^{p}+\sup_{r\in[\frac{1}{2},1)} \frac{1}{2\pi} \int_{0}^{2\pi}|\partial_{r}f(re^{i\theta})|^{p}d\theta\qquad\quad\qquad\qquad\;\;\;\;\qquad\text{(by  \eqref{eq-3.45})}\\
\leq& M_{4}^{p}+\frac{2^{2p-2}}{ \pi}\sup_{r\in[\frac{1}{2},1)} \int_{0}^{2\pi}\left(| \mathscr{J}_{\alpha,1}(re^{i\theta})|^{p}
+  | \mathscr{J}_{\alpha,2}(re^{i\theta})|^{p}\right)d\theta \;\;\text{(by  \eqref{eq-3.48}) }\\
\lesssim & M_{4}^{p}+   \| F \|_{\mathcal{L}^{\infty}(\mathbb{T})}^{p} +  \| \dot F \|_{\mathcal{L}^p(\mathbb{T})}^{p},
\quad\qquad\quad\quad\;\;\text{(by  Lemma \ref{lem-3.3}(1) and (2))}
\end{split}
\end{align*}
 which implies that for $\alpha\in(0,\infty)$ and $p\in[1,\infty)$,
\begin{align}\label{eq-3.49}
\begin{split}
  \|\partial_{r}f \|_{\mathcal{H}_{\mathcal{G}}^{p}(\mathbb{D})}\lesssim  \|F\|_{\mathcal{L}^{\infty}(\mathbb{T})}+\|\dot{F}\|_{\mathcal{L}^{p}(\mathbb{T})}.
\end{split}
\end{align}

\noindent
 (II)  {\bf Estimate on  $\|\partial_{r}f\|_{\mathcal{B}_{\mathcal{G}}^{p}(\mathbb{D})} $}.
 For  $\alpha\in(-1,\infty)$ and $p\in[1,\infty)$,
 obviously, \eqref{eq-3.46} and \eqref{eq-3.48} guarantee that
 \begin{align}\label{eq-3.50}
\begin{split}
B_{p}^{p}\big(\partial_{r}f \big)
\leq& M^{p}_{4}+ \int_{\mathbb{D}\backslash \overline{\mathbb{D}}_{\frac{1}{2}}} |\partial_{r}f(z)|^{p}d\sigma(z) \\
\leq& M^{p}_{4}+2^{2p-1}
\int_{\mathbb{D}\backslash \overline{\mathbb{D}}_{\frac{1}{2}}}\left( |\mathscr{J}_{\alpha,1}(z)|^{p}
+ |\mathscr{J}_{\alpha,2}(z)|^{p}\right)d\sigma(z) .
\end{split}
\end{align}
Notice that   \eqref{eq-3.23} leads to
\begin{align}\label{eq-3.51}
\begin{split}
\int_{\mathbb{D}\backslash \overline{\mathbb{D}}_{\frac{1}{2}}} |\mathscr{J}_{\alpha,1}(z)|^{p}
 d\sigma(z)
\leq \|\mathscr{J}_{\alpha,1}\|^{p}_{   \mathcal{B}_{\mathcal{G}}^{p}(\mathbb{D})}
\leq\|\mathscr{J}_{\alpha,1}\|^{p}_{\mathcal{H}_{\mathcal{G}}^{p}(\mathbb{D})}
\leq|\alpha|^p\| F \|_{\mathcal{L}^{\infty}(\mathbb{T})}^{p}.
\end{split}
\end{align}

Now, by substituting \eqref{eq-3.42} and \eqref{eq-3.51} into \eqref{eq-3.50}, we obtain that for $\alpha=0$ and $p=1$,
\begin{align*}
\begin{split}
B_{p}^{p}\big(\partial_{r}f \big)
\leq& M_{4}^{p}+2^{2p-1}
\int_{\mathbb{D}\backslash \overline{\mathbb{D}}_{\frac{1}{2}}} |\mathscr{J}_{\alpha,2}(z)|^{p}d\sigma(z) \\
\leq& M_{4}^{p}
+ \frac{ \left(\mathcal{C}_{\alpha}\right)^{p}}{\pi^{p}}
\big(2^{4p-2}\log^{p}4+2^{5p-1}\Gamma(1+p)   \big) \|\dot F \|_{\mathcal{L}^p(\mathbb{T})}^{p}.
\end{split}
\end{align*}
Thus for $\alpha=0$ and $p=1$,
\begin{align}\label{eq-3.52}
\begin{split}
  \|\partial_{r}f \|_{\mathcal{B}_{\mathcal{G}}^{p}(\mathbb{D})}\lesssim  \|\dot{F}\|_{\mathcal{L}^{p}(\mathbb{T})}.
\end{split}
\end{align}

Similarly,
by substituting \eqref{eq-3.43} and \eqref{eq-3.51} into \eqref{eq-3.50}, we find that  for $\alpha\in(-1,0)$ and
$p\in\big[1,-\frac{1}{\alpha}\big)$,
\begin{align*}
\begin{split}
 B_{p}^{p}\big(\partial_{r}f \big)
\leq& M_{4}^{p}+2^{2p-1}
\int_{\mathbb{D}\backslash \overline{\mathbb{D}}_{\frac{1}{2}}} \left(|\mathscr{J}_{\alpha,1}(z)|^{p}
+  |\mathscr{J}_{\alpha,2}(z)|^{p} \right)d\sigma(z) \\
\leq& M_{4}^{p}+   2^{2p-1} |\alpha|^{p}\| F \|^{p}_{\mathcal{L}^{\infty}(\mathbb{T})}
+\frac{2^{2 p-1} (\mathcal{C}_{\alpha})^{p}M_{3}}{\pi^{p}}  \big(2^{p} \|\dot F\|^{p}_{\mathcal{L}^{p} (\mathbb{T}) }
+|\alpha|^{p}\|F\|^{p}_{\mathcal{L}^{\infty} (\mathbb{T})}\big),
\end{split}
\end{align*}
 which implies that for $\alpha\in(-1,0)$ and
$p\in\big[1,-\frac{1}{\alpha}\big)$,
 $$
  \|\partial_{r}f \|_{\mathcal{B}_{\mathcal{G}}^{p}(\mathbb{D})}\lesssim \|F\|_{\mathcal{L}^{\infty}(\mathbb{T})}+\|\dot{F}\|_{\mathcal{L}^{p}(\mathbb{T})}.
  $$
This, together with \eqref{eq-3.52}, shows that the second statement of the lemma is true.

\bca\label{case-6} Suppose that $p=\infty$.
\eca

For $\alpha\in(0,\infty)$,
\begin{align*}
\begin{split}
\|\partial_{r}f \|_{\mathcal{H}_{\mathcal{G}}^{\infty}(\mathbb{D})}
=&\|\partial_{r}f \|_{\mathcal{B}_{\mathcal{G}}^{\infty}(\mathbb{D})}
\leq M_{4}+\sup_{z\in \mathbb{D}\backslash \overline{\mathbb{D}}_{\frac{1}{2}}} | \partial_{r}f(z)|
\qquad\qquad\qquad\;\;\text{(by  \eqref{eq-3.47})} \\
 \leq&  M_{4}+2 \sup_{z\in \mathbb{D}\backslash \overline{\mathbb{D}}_{\frac{1}{2}}}\left(|\mathscr{J}_{\alpha,1}(z)|+|\mathscr{J}_{\alpha,2}(z)|\right).
 \qquad\qquad \quad \text{(by  \eqref{eq-3.44})}
\end{split}
\end{align*}
Then it follows from  Lemma \ref{lem-3.3}(1) and (2)  that for $\alpha\in(0,\infty)$,
$$
\|\partial_{r}f \|_{\mathcal{H}_{\mathcal{G}}^{\infty}(\mathbb{D})}\lesssim \|F\|_{\mathcal{L}^{\infty}(\mathbb{T})}+\|\dot{F}\|_{\mathcal{L}^{\infty}(\mathbb{T})} .
$$
This, together with \eqref{eq-3.49}, shows that the first statement of the lemma is true.
The proof of the lemma is complete.
\epf

The following result is based on Lemmas \ref{lem-3.1} and \ref{lem-3.4}.

\begin{lem}\label{lem-3.5}
Under the assumptions of Lemma \ref{lem-3.1}, the  following two statements are true:
\begin{enumerate}
  \item For all $\alpha\in(0,\infty)$ and $p\in[1,\infty]$ or for $\alpha=0$ and all $p \in (1,\infty)$,
then
 $$
  \|\partial_{z} f  \|_{\mathcal{H}^{p}_{\mathcal{G}}(\mathbb{D})}
+ \| \partial_{\overline{z}} f  \|_{\mathcal{H}^{p}_{\mathcal{G}}(\mathbb{D})}
\lesssim  \|  F\|_{\mathcal{L}^{\infty}(\mathbb{T})}+ \|\dot F\|_{\mathcal{L}^p(\mathbb{T})};
$$
  \item For $\alpha=0$ and $p=1$ or for all $\alpha\in(-1,0)$ and
$p\in\big[1,-\frac{1}{\alpha}\big)$,
 then
 $$
 \|\partial_{z} f  \|_{\mathcal{B}^{p}_{\mathcal{G}}(\mathbb{D})}
+ \| \partial_{\overline{z}} f  \|_{\mathcal{B}^{p}_{\mathcal{G}}(\mathbb{D})}
\lesssim   \|  F\|_{\mathcal{L}^{\infty}(\mathbb{T})}+ \|\dot F\|_{\mathcal{L}^p(\mathbb{T})}.
$$
\end{enumerate}
\end{lem}
\bpf For $z=re^{i\theta}\in \mathbb{D}$, since
\begin{align}\label{eq-3.53}
\partial_{\theta}f(z)=i\big(z\partial_{z}f(z)-\overline{z}\partial_{\overline{z}}f(z)\big) \;\;\;\text{and}\;\;\;
\partial_{r}f(z)=\partial_{z}f(z)e^{i\theta}+\partial_{\overline{z}}f(z)e^{-i\theta},
\end{align}
 we see that for any $z \in \mathbb{D}\backslash\{0\}$,
\begin{equation}\label{eq-3.54}
\partial_{z}f(z)=\frac{ r \partial_{r}f(z)-i \partial_{\theta}f(z) }{2z}
\;\; \text{and}\; \;
 \partial_{\overline{z}}f(z)=\frac{ r \partial_{r}f(z) + i  \partial_{\theta}f(z)  }{2\overline{z}}.
\end{equation}
Then we divide the proof into two cases.
\bca\label{case-7} Suppose that $p\in[1,\infty)$.
\eca

For $z=re^{i\theta}\in \mathbb{D} \backslash\{0\}$,
it follows from  \eqref{eq-3.54} and Lemma H that
\begin{align}\label{eq-3.55}
| \partial_{z}f(z)| ^p\leq\frac{1}{2^p}\left(| \partial_{r}f(z)|+  \frac{|\partial_{\theta}f(z)|}{r}  \right)^p
\leq\frac{1}{2}\left(| \partial_{r}f(z)| ^p+\left| \frac{\partial_{\theta}f(z)}{r}\right| ^p\right),
\end{align}
and similarly,
\begin{align}\label{eq-3.56}
| \partial_{\overline{z}}f(z)| ^p\leq\frac{1}{2}\left(| \partial_{r}f(z)| ^p+\left| \frac{\partial_{\theta}f(z)}{r}\right| ^p\right).
\end{align}
Then by \eqref{eq-1.6}, Lemmas \ref{lem-3.1}(2)  and \ref{lem-3.4}, we see that to prove the lemma,
it suffices to show that for any $\alpha\in(-1,\infty)$,
\begin{align*}
\left\| \frac{\partial_{\theta}f }{r}  \right\|_{\mathcal{H}^{p}_{\mathcal{G}}(\mathbb{D})}
=\sup_{r\in[0,1)}\frac{1}{2\pi}\int_{0}^{2\pi}\left| \frac{\partial_{\theta}f(re^{i\theta})}{r}\right |^p\,d\theta
\lesssim  \|\dot F\|_{\mathcal{L}^p(\mathbb{T})},
\end{align*}
where
\begin{align*}
\begin{split}
\lim_{r\rightarrow0^{+}}\left| \frac{\partial_{ \theta }f(re^{i\theta})}{r}\right|
=\left|e^{i\theta} \partial_{z}f(0)- e^{-i\theta} \partial_{\overline{z}}f(0)\right|\leq M_{4}
\end{split}
\end{align*}
and $M_{4}$ is the constant from (\ref{eq-3.47}).

By \eqref{eq-3.53}, we see that
\begin{align}\label{eq-3.57}
\begin{split}
\sup_{r\in[0,\frac{1}{2})}\frac{1}{2\pi}\int_{0}^{2\pi}\left| \frac{\partial_{\theta}f(re^{i\theta})}{r}\right| ^p\,d\theta
\leq&\sup_{r\in[0,\frac{1}{2})}\frac{1}{2\pi}\int_{0}^{2\pi} \| D_{f}(re^{i\theta})\|^p\,d\theta \leq M_{4}^{p} .
\end{split}
\end{align}

Moreover, by Lemma \ref{lem-3.1}(1), we infer that
\begin{align}\label{eq-3.58}
\sup_{r\in[\frac{1}{2},1)}\frac{1}{2\pi}\int_{0}^{2\pi}\left| \frac{\partial_{\theta}f(re^{i\theta})}{r}\right| ^p\,d\theta
\lesssim  \|\dot F\|_{\mathcal{L}^p(\mathbb{T})}^{p}.
\end{align}
The combination  of \eqref{eq-3.57}  and  \eqref{eq-3.58} yields that
 \begin{align}\label{eq-3.59}
 \left\| \frac{\partial_{\theta}f }{r}  \right\|_{\mathcal{B}^{p}_{\mathcal{G}}(\mathbb{D})}
\leq  \left\| \frac{\partial_{\theta}f }{r}  \right\|_{\mathcal{H}^{p}_{\mathcal{G}}(\mathbb{D})}
\lesssim  \|\dot F\|_{\mathcal{L}^p(\mathbb{T})}.
\end{align}
Then  for $\alpha\in(0,\infty)$ and $p\in[1,\infty)$ or for $\alpha=0$ and $p \in (1,\infty)$,
 \begin{align*}
\begin{split}
  \qquad \quad\qquad \|\partial_{z} f&  \|_{\mathcal{H}^{p}_{\mathcal{G}}(\mathbb{D})}
+ \| \partial_{\overline{z}} f  \|_{\mathcal{H}^{p}_{\mathcal{G}}(\mathbb{D})}\\
&\lesssim \left\| \partial_{r}f  \right\|_{\mathcal{H}^{p}_{\mathcal{G}}(\mathbb{D})}+
 \left\| \frac{\partial_{\theta}f }{r}  \right\|_{\mathcal{H}^{p}_{\mathcal{G}}(\mathbb{D})}
  \quad  \;\;\text{(by \eqref{eq-3.55}, \eqref{eq-3.56} and Lemma  H)}\\
&\lesssim   \|  F\|_{\mathcal{L}^{\infty}(\mathbb{T})}+ \|\dot F\|_{\mathcal{L}^p(\mathbb{T})}.
 \quad \; \text{(by \eqref{eq-3.59}, Lemmas \ref{lem-3.1}(2) and \ref{lem-3.4}(1))}
\end{split}
\end{align*}
Similarly, for $\alpha=0$ and $p=1$ or for $\alpha\in(-1,0)$ and
$p\in\big[1,-\frac{1}{\alpha}\big)$,
 \begin{align*}
\begin{split}
  \quad \qquad\qquad \qquad  \|\partial_{z} f & \|_{\mathcal{B}^{p}_{\mathcal{G}}(\mathbb{D})}
+ \| \partial_{\overline{z}} f  \|_{\mathcal{B}^{p}_{\mathcal{G}}(\mathbb{D})}\\
&\lesssim \left\| \partial_{r}f  \right\|_{\mathcal{B}^{p}_{\mathcal{G}}(\mathbb{D})}+
 \left\| \frac{\partial_{\theta}f }{r}  \right\|_{\mathcal{B}^{p}_{\mathcal{G}}(\mathbb{D})}
    \text{(by \eqref{eq-3.55}, \eqref{eq-3.56} and Lemma  H)}\\
&\lesssim  \|  F\|_{\mathcal{L}^{\infty}(\mathbb{T})}+ \|\dot F\|_{\mathcal{L}^p(\mathbb{T})}.
 \quad\qquad\text{(by \eqref{eq-3.59} and  Lemma  \ref{lem-3.4}(2))}
\end{split}
\end{align*}

\bca\label{case-8}
 Suppose that $p=\infty$.
\eca

For $z=re^{i\theta}\in\mathbb{D}\setminus\overline{\mathbb{D}}_{\frac{1}{2}}$ and $\alpha\in(0,\infty)$, we have that
 \begin{align*}
\begin{split}
 \qquad| \partial_{z}f(z)|+| \partial_{\overline{z}}f(z)|
& \leq | \partial_{r}f(z)|+\frac{ | \partial_{\theta}f(z)| }{r} \;\;\;\qquad\qquad\qquad\qquad\qquad\text{(by \eqref{eq-3.54})}\\
&\leq \|\partial_{r}f\|_{\mathcal{H}_{\mathcal{G}}^{\infty}(\mathbb{D})} +2\|\partial_{\theta}f\|_{\mathcal{H}_{\mathcal{G}}^{\infty}(\mathbb{D})}\\
&\lesssim  \|  F\|_{\mathcal{L}^{\infty}(\mathbb{T})}+ \|\dot F\|_{\mathcal{L}^{\infty}(\mathbb{T})}. \;\;\;\;\text{(by Lemmas \ref{lem-3.1}(1) and \ref{lem-3.4}(1))}
\end{split}
\end{align*}
Since $\partial_{z}f\in \mathcal{C}(\overline{\mathbb{D}}_{\frac{1}{2}})$ and $\partial_{\overline{z}}f\in \mathcal{C}(\overline{\mathbb{D}}_{\frac{1}{2}})$,
we get that
 \begin{align*}
\begin{split}
\| \partial_{z}f \|_{\mathcal{H}_{\mathcal{G}}^{\infty}}+\| \partial_{\overline{z}}f \|_{\mathcal{H}_{\mathcal{G}}^{\infty}}
\lesssim  \|  F\|_{\mathcal{L}^{\infty}(\mathbb{T})}+ \|\dot F\|_{\mathcal{L}^{\infty}(\mathbb{T})}.
\end{split}
\end{align*}
Hence the proof of this lemma is complete.
\epf

Based on Lemma \ref{lem-3.1}, we are going to establish the following result concerning $(K,K')$-elliptic mappings.

\begin{lem}\label{lem-3.6}
Under the assumptions of Lemma \ref{lem-3.1}, further, if $f$ is a $(K,K')$-elliptic mapping in $\mathbb{D}$,
then for all $\alpha\in(-1,0)$ and $p\in[1,\infty]$ or for $\alpha=0$ and $p\in\{1,\infty\}$,
 $$\|\partial_{r}f\|_{\mathcal{H}_{\mathcal{G}}^{p}( \mathbb{D} )}
+\|\partial_{z}f\|_{\mathcal{H}_{\mathcal{G}}^{p}( \mathbb{D} )}
+
\|\partial_{\overline{z}}f\|_{\mathcal{H}_{\mathcal{G}}^{p}( \mathbb{D} )}
\lesssim \| \dot F\|_{\mathcal{L}^p(\mathbb{T})}.
$$
\end{lem}
\bpf In order to prove the lemma, we divide the proof into two cases.

\bca\label{case-9} Suppose that $p\in [1,\infty).$
\eca

Since $f$ is a $(K,K')$-elliptic mapping in $\mathbb{D}$, we obtain from \cite[Inequality (2.7)]{csw}  that
\begin{equation}\label{eq-3.60}
l^p\big(D_{f}(z)\big)\geq\frac{1}{2^{p-1}K^p}\| D_{f}(z)\| ^p-\left(\frac{\sqrt{K' }}{K } \right)^{p}.
\end{equation}
For $z=re^{i\theta}\in\mathbb{D}$, by calculations, we get that
\begin{align*}
\begin{split}
\qquad\qquad\|\dot F\|_{\mathcal{L}^p(\mathbb{T})}^{p}\geq& \frac{1}{2\pi} \int_{0}^{2\pi} | \partial_{\theta}f(re^{i\theta})| ^p\,d\theta\;\;\;\;\;\;\;\;\;\;\;\;\;\;\;\;\;\;\qquad\qquad \quad\qquad\quad\text{(by \eqref{eq-3.7})}\\
\geq& \frac{r^p}{2\pi} \int_{0}^{2\pi}l^p(D_{f}\big(re^{i\theta})\big)\,d\theta\quad\;
\qquad\qquad\qquad\;\; \text{(by \eqref{eq-1.1} and (\ref{eq-3.53}))}\\
\geq&\frac{r^p}{ 2^{p  }K^p  \pi}\int_{0}^{2\pi}\| D_{f}(re^{i\theta})\| ^p\,d\theta
- \left(\frac{\sqrt{K' }}{K } \right)^{p},\;\;\qquad\qquad\text{(by (\ref{eq-3.60}))}
\end{split}
\end{align*}
which means that
\begin{equation}\label{eq-3.61}
 \sup\limits_{r\in[\frac{1}{2},1)} \left(\frac{1}{2\pi} \int_{0}^{2\pi}\| D_{f}(re^{i\theta})\| ^p\,d\theta \right)^{\frac{1}{p}}
\leq \left( 2^{2p-1}K^p \| \dot F\|_{\mathcal{L}^p(\mathbb{T})}^p+ 2^{2p-1} K'^{\frac{p}{2}}\right)^{\frac{1}{p}} .
\end{equation}
Recall that $\| D_{f}\|=|f_{z}|+|f_{\overline{z}}|$ is continuous in $\overline{\mathbb{D}}_{\frac{1}{2}}$.
Therefore,
 \begin{equation}\label{eq-3.62}
\sup\limits_{r\in [0,\frac{1}{2}]}\left(\frac{1}{2\pi} \int_{0}^{2\pi}\| D_{f}(re^{i\theta})\| ^p\,d\theta \right)^{\frac{1}{p}}
\leq
\,M_{4},
\end{equation} where $M_{4}$ is the constant from (\ref{eq-3.47}).
Combining \eqref{eq-3.61} and \eqref{eq-3.62}, we see that
 $$
\|D_{f}\|_{\mathcal{H}_{\mathcal{G}}^{p}( \mathbb{D} )}\lesssim \| \dot F\|_{\mathcal{L}^p(\mathbb{T})}.
$$
Then \eqref{eq-1.1} yields that
 \begin{equation}\label{eq-3.63}
\|\partial_{z}f\|_{\mathcal{H}_{\mathcal{G}}^{p}( \mathbb{D} )}
+
\|\partial_{\overline{z}}f\|_{\mathcal{H}_{\mathcal{G}}^{p}( \mathbb{D} )}
\leq 2\|D_{f}\|_{\mathcal{H}_{\mathcal{G}}^{p}( \mathbb{D} )}
\lesssim \| \dot F\|_{\mathcal{L}^p(\mathbb{T})}.
\end{equation}
 Further, for $z=re^{i\theta}\in \mathbb{D}$,
it follows from \eqref{eq-3.53}, \eqref{eq-3.63} and Lemma H that
 \begin{align*}
\begin{split}
 \|\partial_{r}f\|_{\mathcal{H}_{\mathcal{G}}^{p}( \mathbb{D} )}
&\leq  2^{\frac{p-1}{p}}\sup\limits_{r\in[0,1)} \left(\frac{1}{2\pi} \int_{0}^{2\pi}\big(|\partial_{z}f(re^{i\theta})|^p+|\partial_{\overline{z}}f(re^{i\theta})|^p\big)\,d\theta \right)^{\frac{1}{p}}\\
&\leq 2^{\frac{p-1}{p}}\big(\|\partial_{z}f\|_{\mathcal{H}_{\mathcal{G}}^{p}( \mathbb{D} )}
+
\|\partial_{\overline{z}}f\|_{\mathcal{H}_{\mathcal{G}}^{p}( \mathbb{D} )}\big)
 \lesssim \| \dot F\|_{\mathcal{L}^p(\mathbb{T})},
\end{split}
\end{align*}
which is what we need.

\bca\label{case-10} Suppose that $p=\infty$.
\eca
For $z=re^{i\theta}\in \mathbb{D}$, we obtain from \eqref{eq-1.1},  \eqref{eq-3.53}  and \cite[Lemma 1.1]{clsw} that
$$
H_{\infty}(\partial_{\theta}f )
\geq | \partial_{\theta}f(z)|
 \geq |z| \cdot  l\big(D_{f}(z)\big)
 \geq \frac{|z|}{K}\left(\| D_{f}(z)\|-\sqrt{K'}\right).
$$
Then Lemma \ref{lem-3.1}(1) leads to
$$
\sup_{z\in \mathbb{D}} \left(| z| \cdot\| D_{f}(z)\|\right)
\leq \sqrt{K'}+K\cdot H_{\infty}(\partial_{\theta}f )
\lesssim \| \dot F\|_{\mathcal{L}^{\infty}(\mathbb{T})}.
$$
Therefore, for $z\in \mathbb{D}\setminus\overline{\mathbb{D}}_{\frac{1}{2}}$,
we infer from \eqref{eq-1.1} and \eqref{eq-3.53} that
\begin{align*}
\begin{split}
|\partial_{r}f(z)|
&\leq| \partial_{z}f(z)| +|\partial_{\overline{z}}f(z)| \leq 2\sup\limits_{z\in \mathbb{D}\setminus\overline{\mathbb{D}}_{\frac{1}{2}}}(|z|\cdot \| D_{f}(z)\|) \lesssim \| \dot F\|_{\mathcal{L}^{\infty}(\mathbb{T})}.
\end{split}
\end{align*}
Since
 $$
  \max\{|\partial_{r}f(z)|+|\partial_{z}f(z)|+|\partial_{\overline{z}}f(z)|:z\in\overline{\mathbb{D}}_{\frac{1}{2}} \}
  \leq 2M_{4},
 $$
 we see that
$$
\|\partial_{r}f\|_{\mathcal{H}_{\mathcal{G}}^{\infty}(\mathbb{D})}
+\|\partial_{z}f\|_{\mathcal{H}_{\mathcal{G}}^{\infty}(\mathbb{D})}
+\|\partial_{\overline{z}}f \|_{\mathcal{H}_{\mathcal{G}}^{\infty}(\mathbb{D})}
\lesssim  \| \dot F\|_{\mathcal{L}^{\infty}(\mathbb{T})},
$$
 where $M_{4}$ is the constant from (\ref{eq-3.47}). The proof of the lemma is complete.
\epf

\subsection{Proof of Theorem \ref{thm-1.1}}
Theorem \ref{thm-1.1}(1) (resp. Theorem \ref{thm-1.1}(2); Theorem \ref{thm-1.1}(3)) follows from Lemma \ref{lem-3.1}(1)
(resp. Lemmas \ref{lem-3.1}(2), \ref{lem-3.4}(1) and \ref{lem-3.5}(1); Lemma \ref{lem-3.6}).
\qed

\subsection{Proof of Theorem \ref{thm-1.3}}
Theorem 1.3 follows from Lemmas \ref{lem-3.4}(2) and \ref{lem-3.5}(2).
\qed


\section{Proofs  of Theorems \ref{thm-1.2} and \ref{thm-1.4} }\label{sec-4}
The purpose of this section is to prove Theorems \ref{thm-1.2} and \ref{thm-1.4}.
 The proofs are based on the following examples.


\subsection{Examples}

\begin{eg}\label{ex-4.1}
Suppose that $\alpha\in(-1,0)$ and
 $n\in \mathbb{N } $. Let
\begin{equation*}
f(z)
= \;_{2}F_{1}\left(-\frac{\alpha}{2}, n-\frac{\alpha}{2}; n+1; | z|^2\right)z^n
\end{equation*}
in $ \overline{\mathbb{D}}$, and let
\begin{equation*}
F(\xi)=\;_{2}F_{1}\left(-\frac{\alpha}{2}, n-\frac{\alpha}{2}; n+1; 1\right)\xi^n
\end{equation*}
on $\mathbb{T}$.
Then the following five statements are true:
 \begin{enumerate}
 \item\label{ex-4.1-1}
   $F$ is an absolutely continuous function on $\mathbb{T}$, $\dot F \in \mathcal{C}^{\infty}(\mathbb{T})$  and  $f_{r}  \rightarrow F $
in $\mathcal{D}'(\mathbb{T})$ as $r\rightarrow 1^{-}$;
  \item\label{ex-4.1-2}
    $f \in \mathcal{C}^{\infty}(\mathbb{D})$,  $T_{\alpha}f =0$ and $f =\mathcal{K}_{\alpha}[F ]$ in $\mathbb{D}$;
 \item\label{ex-4.1-3}
  $f $ is not $(K,K')$-elliptic for any $K\geq1$ and $K'\geq 0$;
  \item \label{ex-4.1-4}
    For any $p\in(0, \infty]$,
$$
 H_p(\partial_{r}f)=H_p(\partial_{z}f)=H_p(\partial_{\overline{z}}f)=\infty;
$$
 \item \label{ex-4.1-5}
 For any $p\in[-\frac{1}{\alpha}, \infty]$,
$$
B_p( \partial_{r}f)=B_p( \partial_{z}f)=B_p( \partial_{\overline{z}}f)=\infty.
$$
\end{enumerate}
\end{eg}

\bpf
We firstly prove $(1)$. Obviously, $F $ is an absolutely continuous function on $\mathbb{T}$ and $\dot F \in \mathcal{C}^{\infty}(\mathbb{T})$. To check the remaining assertion in $(1)$,
we first show that $f\in \mathcal{C} (\overline{\mathbb{D}})$.

Observe that
\begin{equation}\label{eq-4.1}
f(z)= \sum_{k=0}^{\infty}a_{n,k}|z|^{2k}z^{n},
\end{equation}
where
$$
a_{n,k}=\frac{(-\frac{\alpha}{2})_{k} (n-\frac{\alpha}{2})_{k}}{(n+1)_{k}k!}.
$$
Since $a_{n,k}>0$ for each $k\in \mathbb{N}_{0}$  and
$$
\sum_{k=0}^{\infty}a_{n,k}=\;_{2}F_{1}\left(-\frac{\alpha}{2}, n-\frac{\alpha}{2}; n+1; 1\right)<\infty,
$$
 we see that
the series $\;_{2}F_{1}\left(-\frac{\alpha}{2}, n-\frac{\alpha}{2}; n+1; s\right)$ converges absolutely and uniformly in $ [0,1]$,
which implies that
$ \;_{2}F_{1}\left(-\frac{\alpha}{2}, n-\frac{\alpha}{2}; n+1; s\right)$ is continuous in $[0,1]$.
Therefore, $f\in \mathcal{C} (\overline{\mathbb{D}})$.
Obviously, for any $\xi\in \mathbb{T}$,
$$
\lim_{r\rightarrow1^{-}} f(r\xi)=F(\xi).
$$

Now, it follows from \eqref{eq-2.4} that for every $\varphi\in \mathcal{C}^{\infty}(\mathbb{T})$,
$$
\lim_{r\rightarrow1^{-}} \left\langle  f_{r} ,\varphi \right\rangle
=\lim_{r\rightarrow1^{-}}\frac{1}{2\pi}\int_{\mathbb{T}} \varphi(e^{i\theta})f(re^{i\theta})d\theta
=\frac{1}{2\pi} \int_{\mathbb{T}} \varphi(e^{i\theta})F( e^{i\theta})d\theta
= \left\langle F ,\varphi  \right\rangle,
$$
i.e., $ f_{r}\rightarrow F $
in $\mathcal{D}'(\mathbb{T})$ as $r\rightarrow 1^{-}$.
This proves $(1)$ in the example.

Secondly, we prove $(2)$.
Since elementary calculations show that the convergent radius of the series in \eqref{eq-4.1} is $1$,
and so, $f\in  \mathcal{C}^{\infty}(\mathbb{D})$. Then we know from \cite[Theorem $2.1$]{ao} that $T_{\alpha}(f)=0$ in $\mathbb{D}$.
Because we have known that
$f\in  \mathcal{C}^{\infty}(\mathbb{D})$,  $T_{\alpha}(f)=0$ in $\mathbb{D}$ and
$ f_{r}\rightarrow F $
in $\mathcal{D}'(\mathbb{T})$ as $r\rightarrow 1^{-}$,
we see from  \cite[Theorem 3.3]{ao}   that $f=\mathcal{K}_{\alpha}[F]$ in $ \mathbb{D} $.
 Hence the statement $(2)$ in the example is true.

Thirdly, we prove $(3)$. Elementary computations show that
\begin{align}\label{eq-4.2}
\begin{split}
 \partial_{z}f(z)
=&\frac{ \alpha(\alpha-2n)}{4(n+1)}
r^{n+1} \xi^{n-1} \mathscr{E}_{1}(r)
 +n r^{n-1}\xi^{n-1}\mathscr{E}_{2}(r),
\end{split}
\end{align}
\begin{equation}\label{eq-4.3}
 \partial_{\overline{z}}f(z)
=\frac{ \alpha(\alpha-2n)}{4(n+1)}
r^{n+1}\xi^{n+1}\mathscr{E}_{1}(r)
\end{equation}
and
\begin{align}\label{eq-4.4}
\begin{split}
\partial_{r}f(z)
= &\frac{ \alpha(\alpha-2n)}{2(n+1)}r^{n+1} \xi^{n}\mathscr{E}_{1}(r)
 +nr^{n-1}\xi^{n}\mathscr{E}_{2}(r),
\end{split}
\end{align}
where $z=r\xi$ with $r\in[0,1)$,
\begin{align}\label{eq-4.5}
\begin{split}
\mathscr{E}_{1}(r)=\;_{2}F_{1}\left(1-\frac{\alpha}{2}, n+1-\frac{\alpha}{2}; n+2; r^2\right)
\end{split}
\end{align}
and
$$\mathscr{E}_{2}(r)=\;_{2}F_{1}\left(-\frac{\alpha}{2}, n-\frac{\alpha}{2}; n+1; r^2\right).$$

For the quantity $\mathscr{E}_{1}(r)$, it follows from \eqref{eq-2.2} and  \eqref{eq-2.3} that
\begin{align}
\begin{split}\label{eq-4.6}
  \lim_{r\rightarrow1^{-}}\mathscr{E}_{1}(r)
&= \lim_{r\rightarrow1^{-}}(1-r^{2})^{\alpha} \;_{2}F_{1}\left(n+1+\frac{\alpha}{2}, 1+\frac{\alpha}{2}; n+2, r^2\right)\\
&= \frac{\Gamma(n+2)\Gamma(-\alpha)}{\Gamma(1-\frac{\alpha}{2})\Gamma(n+1-\frac{\alpha}{2})  } \lim_{r\rightarrow1^{-}}(1-r^{2})^{\alpha}
=\infty,
\end{split}
\end{align}
and for the quantity $\mathscr{E}_{2}(r)$, we infer from \eqref{eq-2.1} and \eqref{eq-2.2} that
\begin{align}
\begin{split}\label{eq-4.7}
0<&\mathscr{E}_{2}(r)
\leq \; _{2}F_{1}\left(-\frac{\alpha}{2}, n-\frac{\alpha}{2}; n+1; 1\right)
 =  \frac{\Gamma(n+1)\Gamma(1+\alpha)}
{\Gamma(n+1+\frac{\alpha}{2})\Gamma(1+\frac{\alpha}{2}) }
<\infty.
\end{split}
\end{align}

Now, by (\ref{eq-1.1}), \eqref{eq-4.2} and \eqref{eq-4.3}, we have that
\begin{align}
\begin{split}\label{eq-4.8}
\|D_{f} (z)\|
= | \partial_{z}f(z)|+| \partial_{\overline{z}}f(z)|
=\frac{ \alpha(\alpha-2n)}{2(n+1)}
r^{n+1}\mathscr{E}_{1}(r)
 +nr^{n-1}\mathscr{E}_{2}(r)
\end{split}
\end{align}
and
\begin{align*}
\begin{split}
 l\big(D_{f} (z)\big)
= | \partial_{z}f(z)|-| \partial_{\overline{z}}f(z)|
= nr^{n-1}\mathscr{E}_{2}(r).
\end{split}
\end{align*}
Then for any $K\geq1$,
\begin{align}
\begin{split}\label{eq-4.9}
 \|D_{f} (z)\|-Kl\big(D_{f} (z)\big) 
 =\frac{ \alpha(\alpha-2n)}{2(n+1)}
 r^{n+1}\mathscr{E}_{1}(r)
  +n(1-K) r^{n-1}\mathscr{E}_{2}(r).
\end{split}
\end{align}
Substituting \eqref{eq-4.6} and \eqref{eq-4.7} into \eqref{eq-4.8} and \eqref{eq-4.9}, respectively, we see that
\begin{align*}
\begin{split}
\lim_{|z|\rightarrow1^{-}}\left(\|D_{f}(z)\|^2-KJ_{f}(z)  \right)
&=\lim_{|z|\rightarrow1^{-}}\|D_{f}(z)\|\left(\|D_{f}(z)\|-Kl\big(D_{f}(z)\big)\right)
=\infty.
\end{split}
\end{align*}
This shows that $f $ is not $(K,K')$-elliptic for any $K\geq 1$ and $K'\geq 0$, and thus, the third statement in the example is proved.

Fourthly, we prove $(4)$. Since for $i\in \{1,2\}$ and for any $r\in [0,1)$, $\mathscr{E}_{i}(r)\geq 0$, we see that for any  $p\in(0,\infty)$  and $r\in[0,1)$,
\eqref{eq-4.4} and \eqref{eq-4.7} guarantee that
\begin{align*}
\frac{1}{2\pi}\int_{0}^{2\pi}| \partial_{r}f(re^{i\theta})| ^p\,d\theta
=&\frac{1}{2\pi}\int_{0}^{2\pi}\left|\frac{\alpha(\alpha -2n )r^{n+1}}{2(n+1)} \mathscr{E}_{1}(r) +nr^{n-1}\mathscr{E}_{2}(r)\right| ^p\,d\theta\\
\geq & \left(\frac{\alpha(\alpha-2n)r^{n+1}}{2(n+1)} \mathscr{E}_{1}(r) \right)^p
\geq0.
\end{align*}
Then  \eqref{eq-4.6} yields  that for any $p\in(0,\infty)$,
\begin{align*}
H_p( \partial_{r}f)
= \sup_{r\in[0,1)}\frac{1}{2\pi}\int_{0}^{2\pi}| \partial_{r}f(re^{i\theta})|^p\,d\theta
 \geq & \lim_{r\rightarrow1^{-}}\left( \frac{\alpha(\alpha-2n) }{2(n+1)}\mathscr{E}_{1}(r) \right)^p
 = \infty.
\end{align*}
Thus
 the fact
 $$ H_{\infty}(\partial_{r}f)= \sup_{z\in \mathbb{D}}|\partial_{r}f(z)| \geq H_1(\partial_{r}f)$$
guarantees that for any $p\in(0,\infty]$,
\begin{align}
\begin{split}\label{eq-4.10}
 H_p(\partial_{r}f)=\infty.
\end{split}
\end{align}

Moreover, we deduce from \eqref{eq-4.2},  \eqref{eq-4.3}
and similar arguments as in the proof of \eqref{eq-4.10} that for any $p\in(0,\infty]$,
$$
 H_p(\partial_{z}f)=H_p(\partial_{\overline{z}}f)=\infty.
$$
These complete the proof of  the fourth statement in the example.

Fifthly, we prove $(5)$. For $p\in[-\frac{1}{\alpha},\infty)$, by \eqref{eq-4.3},
 we find that
\begin{align*}
\int_{\mathbb{D}} |\partial_{\overline{z}}f(z)|^p\,d\sigma(z)
=&2\int_{0}^{1}\left( \frac{\alpha(\alpha -2n )r^{n+1}}{4(n+1)} \mathscr{E}_{1}(r) \right)^p r \,dr\\
=& \frac{2 (\alpha^{2}-2n\alpha)^{p} }{4^{p}(n+1)^p } \int_{0}^{1}r^{(n+1)p+1}  \mathscr{E}_{1}^{p}(r) \,dr,
\end{align*} where $z=re^{i\theta}$.
Since \eqref{eq-2.3} and \eqref{eq-4.5} implies that
$$
\mathscr{E}_{1}(r)
 =(1-r^2)^{\alpha  }
  \;_{2}F_{1}\left(n+1+\frac{\alpha}{2}, 1+\frac{\alpha}{2}; n+2; r^2 \right),
$$
  we have
\begin{align*}
\int_{\mathbb{D}}&|\partial_{\overline{z}}f(z)|^p\,d\sigma(z)\\
=&  \frac{2 (\alpha^{2}-2n\alpha)^{p} }{4^{p}(n+1)^p } \int_{0}^{1}r^{(n+1)p+1}(1-r^2)^{\alpha p}
\left( \;_{2}F_{1}\left(n+1+\frac{\alpha}{2}, 1+\frac{\alpha}{2}; n+2; r^2 \right)\right)^{p} \,dr.
\end{align*}
Then the fact $_{2}F_{1}\left(n+1+\frac{\alpha}{2}, 1+\frac{\alpha}{2}; n+2; r^2 \right)\geq 1$ in $[0,1)$ leads to
\begin{align*}
B_p( \partial_{\overline{z}}f)=\int_{\mathbb{D}} |\partial_{\overline{z}}f(z)|^p\,d\sigma(z)
\geq    \frac{2 (\alpha^{2}-2n\alpha)^{p} }{4^{p}(n+1)^p } \int_{0}^{1}r^{(n+1)p+1}(1-r^2)^{\alpha p} \,dr.
\end{align*}
But for $\alpha\in(-1,0)$ and $p\in[-\frac{1}{\alpha},\infty)$,
the integral $\int_{0}^{1}r^{(n+1)p+1}(1-r^2)^{\alpha p} \,dr$ diverges, and thus,
$$
B_p( \partial_{\overline{z}}f)=\infty.
$$
 This, together with the fact $B_{\infty}( \partial_{\overline{z}}f)\geq B_p( \partial_{\overline{z}}f)$,
yields that for any $p\in[-\frac{1}{\alpha},\infty]$,
$$B_p( \partial_{\overline{z}}f)=\infty.$$

Similarly, for any
 $\alpha\in(-1,0)$
 and $p\in[-\frac{1}{\alpha},\infty]$, also
we have
$$
B_p( \partial_{r}f)=\infty
\;\;\text{and}\;\;
B_p( \partial_{z}f)=\infty.
$$
These demonstrate that the fifth statement in the example is true.
\epf

\begin{eg}\label{ex-4.2}
For $\theta\in[-\pi,\pi]$, set
$$
\phi(\theta)=\begin{cases}
\displaystyle 1+ \theta(\pi+1)/\pi,
 \;\;&-\pi\leq \theta<0, \\
\displaystyle 1+\theta(\pi-1)/\pi, \;\;&0\leq \theta\leq \pi.
\end{cases}
$$
Let $$F  (e^{i\theta})=e^{i\phi(\theta)}$$ on $\mathbb{T}$, and let $$f =\mathcal{K}_{0}[F]$$ in $\mathbb{D}$.
 Then the following five statements are true:
\begin{enumerate}
\item\label{ex-4.2-1}  $F $ is absolutely continuous on $\mathbb{T}$ and $\dot{F}  \in \mathcal{L}^{\infty}(\mathbb{T})$;
\item\label{ex-4.2-2}
$f$ is harmonic in $\mathbb{D}$ and $f_{r} \rightarrow F $
in $\mathcal{D}'(\mathbb{T})$ as $r\rightarrow 1^{-}$;
\item\label{ex-4.2-3}
$ H_{\infty}(\partial_{r}f)=H_{\infty}(\partial_{z}f)=H_{\infty}(\partial_{\overline{z}}f) =\infty$;
\item\label{ex-4.2-4}
$ B_{\infty}(\partial_{r}f)=B_{\infty}(\partial_{z}f)=B_{\infty}(\partial_{\overline{z}}f) =\infty$;
\item\label{ex-4.2-5}
 $f $ is not $(K,K')$-elliptic for any $K\geq1$ and $K'\geq 0$.
\end{enumerate}
\end{eg}

\bpf
Obviously, $F$ is  absolutely continuous on $\mathbb{T}$,
 $\dot{F} \in \mathcal{L}^{\infty}(\mathbb{T})$
and $f$ is harmonic in $\mathbb{D}$.
By the proof of  \cite[Example 4.1]{zjf},
we know that
$$ B_{\infty}(\partial_{r}f)=
B_{\infty}(\partial_{z}f)= B_{\infty}(\partial_{\overline{z}}f)=\infty,$$
which yields that
$$ H_{\infty}(\partial_{r}f)=
H_{\infty}(\partial_{z}f)=H_{\infty}(\partial_{\overline{z}}f)=\infty.$$

 Since
$$
\text{esssup}_{z\in \mathbb{D}}\left|\frac{f_{\overline{z}}(z)}{f_{z}(z)}\right|=1
$$
 (cf. \cite[Example 4.1]{zjf}),
$B_{\infty}(\partial_{z}f)=B_{\infty}(\partial_{\overline{z}}f)=\infty$,
  we know that $f$ is not $(K,K')$-elliptic for any $K\geq1$ and $K'\geq 0$.

As  \cite[Theorem 3.1.6]{pav} ensures that
 $$
 \lim_{z\rightarrow\xi\in\mathbb{ T}}f(z)=\lim_{z\rightarrow\xi}\mathcal{K}_{0}[F](z)=F(\xi),
 $$
 similar arguments as in the proof of Example \ref{ex-4.1}(1) guarantee that
$f_{r} \rightarrow F $
in $\mathcal{D}'(\mathbb{T})$ as $r\rightarrow 1^{-}$. These complete the proof of the example.
\epf

\begin{eg}\label{ex-4.3}
Let
\begin{align}\label{eq-4.11}
\begin{split}
 f (z)=\text{Im} \left(\sum_{n=2}^{\infty} \frac{z^{n}}{n\log n} \right)
\end{split}
\end{align}
  in $\mathbb{D}$,
 and let
\begin{align*}
\begin{split}
F(e^{i\theta})= \sum_{n=1}^{\infty} \frac{\sin (n\theta)}{n\log n}
\end{split}
\end{align*}
on $\mathbb{T}$.
Then the following four statements are true:
\begin{enumerate}
\item \label{ex-4.3-1}
$F $ is absolutely continuous on $\mathbb{T}$ and $\dot{F}  \in \mathcal{L}^{1}(\mathbb{T})$;
\item \label{ex-4.3-2}
$f$ is harmonic in $\mathbb{D}$, $f =\mathcal{K}_{0}[F ]$ in $\mathbb{D}$  and $f_{r} \rightarrow F $
in $\mathcal{D}'(\mathbb{T})$ as $r\rightarrow 1^{-}$;
\item\label{ex-4.3-3}
$ H_{1}(\partial_{r}f)=H_{1}(\partial_{z}f)=H_{1}(\partial_{\overline{z}}f)=\infty $;
\item\label{ex-4.3-4}
 $f $ is not a $(K,K')$-elliptic mapping for any $K\geq1$ and $K'\geq 0$.
\end{enumerate}
\end{eg}

\bpf
(1) Notice that the series $\sum_{n=1}^{\infty} \frac{\sin (n\theta)}{n\log n}$   converges uniformly in $\mathbb{R}$.
  Hence $F \in \mathcal{C}(\mathbb{T})$.
Since $ \sum_{n=1}^{\infty} \frac{\cos(n \theta )}{ \log n}$   converges uniformly in any $[a,b]\subset (0,2\pi)$,
we get
\begin{align*}
\begin{split}
\dot{F}  (e^{i\theta})
=\frac{d}{d \theta}\sum_{n=1}^{\infty} \frac{\sin(n\theta) }{n \log n}
=\sum_{n=1}^{\infty} \frac{\cos(n\theta) }{ \log n}
\end{split}
\end{align*}
in $(0,2\pi)$.
By \cite[p.63]{dp}, we know that $\dot{F} \in \mathcal{L}^{1}(\mathbb{T})$.

 Further, since $F \in \mathcal{C}(\mathbb{T})$,
  $\dot{F}(e^{i\theta})$ exists for every $\theta\in(0,2\pi)$
  and $\dot{F} \in \mathcal{L}^{1}(\mathbb{T})$, we know from \cite[p.153]{ss} that
$F(e^{i \theta})$ is absolutely continuous in $[0,2\pi]$.
This yields that $F $ is absolutely continuous on $\mathbb{T}$.
The first statement in the example is proved.

(2) Elementary calculations show that both the convergent radii of the series  $\sum_{n=2}^{\infty} \frac{z^{n}}{n\log n}$ and
 $\sum_{n=2}^{\infty} \frac{z^{n-1}}{\log n}$ are equal to $1$,
which yields that both series $\sum_{n=2}^{\infty} \frac{z^{n}}{n\log n}$ and
 $\sum_{n=2}^{\infty} \frac{z^{n-1}}{\log n}$ are uniformly convergent in
 $\overline{\mathbb{D}}_{\rho_{0}}$
for any $\rho_{0}\in(0,1)$.
 Then $\sum_{n=2}^{\infty} \frac{z^{n}}{n\log n}$ and
 $\sum_{n=2}^{\infty} \frac{z^{n-1}}{\log n}$ are analytic in $  \mathbb{D} $, and
  $$
  \frac{d}{dz}\left(\sum_{n=2}^{\infty} \frac{z^{n}}{n\log n}\right)=\sum_{n=2}^{\infty} \frac{z^{n-1}}{\log n}.
  $$
 Therefore, $f$ defined in \eqref{eq-4.11} is  harmonic in $  \mathbb{D} $.

Since $F $ is absolutely continuous on $\mathbb{T}$ and
 $\lim_{r\rightarrow1^{-}}f(r\xi)= F(\xi)$
  for any $\xi\in  \mathbb{T}$,
similar arguments as   in the proof of Example \ref{ex-4.1}(1)
show that
$f_{r} \rightarrow F$
in $\mathcal{D}'(\mathbb{T})$ as $r\rightarrow 1^{-}$.

Because we have known that
$f$ is harmonic in $\mathbb{D}$
and
$f_{r} \rightarrow F$
in $\mathcal{D}'(\mathbb{T})$ as $r\rightarrow 1^{-}$,
we see from \cite[Theorem 3.3]{ao}   that $f=\mathcal{K}_{0}[F ]$ in $\mathbb{D} $.
The second statement in the example is true.

(3) For $z=re^{i\theta}\in \mathbb{D}$, by \eqref{eq-4.11}, we find that
 $$
 \partial_{\theta}f(z)= \sum_{n=2}^{\infty} \frac{r^{n}\cos(n\theta)}{ \log n}
 \;\;
 \text{and}\;\;
 \partial_{r}f(z)= \sum_{n=2}^{\infty} \frac{r^{n-1}\sin(n\theta)}{ \log n}.
 $$
Since $f$ satisfies the assumptions in Lemma \ref{lem-3.1}, we see from   Lemma \ref{lem-3.1}(1)  that $\partial_{\theta}f \in \mathcal{H}_{\mathcal{G}}^{1}(\mathbb{\mathbb{D}})$.
By  \cite[p.63]{dp}, we know that
$$
\partial_{\theta}f(z)+ir \partial_{r}f(z)
=\sum_{n=2}^{\infty} \frac{z^{n}}{\log n}
\not\in \mathcal{H}_{\mathcal{G}}^{1}(\mathbb{D}),
$$
which implies that
$r \partial_{r}f(z)
\not\in \mathcal{H}_{\mathcal{G}}^{1}(\mathbb{D})$.
Further,
\begin{align}\label{eq-4.12}
\begin{split}
 H_{1}(\partial_{r}f)
=\sup_{r\in[0,1)}\int_{0}^{2\pi}|   \partial_{r}f(re^{i\theta})| d\theta
\geq
 \sup_{r\in[0,1)}\int_{0}^{2\pi}r| \partial_{r}f(re^{i\theta})| d\theta
=\infty.
\end{split}
\end{align}

Next, we are going to prove
$  H_{1}(\partial_{z}f)= H_{1}(\partial_{\overline{z}}f) =\infty$.

By the obvious fact $\partial_{r}f\in \mathcal{C}(\overline{\mathbb{D}}_{\frac{1}{2}})$ and \eqref{eq-4.12},
 we obtain that
$$
  \sup_{\frac{1}{2}\leq r<1}\int_{0}^{2\pi} \big|\partial_{r}f(re^{i\theta})\big|\,d\theta= \infty.
$$
Since Lemma \ref{lem-3.1}(1) ensures that
$\|\partial_{\theta}f\|_{\mathcal{H}_{\mathcal{G}}^{1}(\mathbb{D})}
\lesssim \|\dot{F}\|_{\mathcal{L}^{1}(\mathbb{T})}$,
we know from  \eqref{eq-3.54}   that
\begin{align*}
\begin{split}
\sup_{\frac{1}{2}\leq r<1} \int_{0}^{2\pi} \big|\partial_{z}f(re^{i\theta})\big| \,d\theta
\geq& \frac{1}{2} \sup_{\frac{1}{2}\leq r<1}
\int_{0}^{2\pi}\left(\big|\partial_{r}f(re^{i\theta})\big|-2\big|\partial_{\theta}f(re^{i\theta})\big|\right)\,d\theta
= \infty,
\end{split}
\end{align*}
which gives
\begin{align}\label{eq-4.13}
\begin{split}
H_{1}(\partial_{z}f) =\infty.
\end{split}
\end{align}

For any $z\in \mathbb{D}$, we infer from \eqref{eq-4.11}  that
$$
 \partial_{z}f(z)=\frac{1}{2i}  \sum_{n=2}^{\infty} \frac{z^{n-1}}{ \log n}
\;\;\text{and}\;\;
 \partial_{\overline{z}}f(z)=\frac{1}{2i}  \sum_{n=2}^{\infty} \frac{\overline{z}^{n-1}}{ \log n},
$$
which ensures that
\begin{align}\label{eq-4.14}
\begin{split}
|\partial_{z}f(z)|=|\partial_{\overline{z}}f(z)|
\quad\text{and}\quad
J_{f}(z)\equiv0.
\end{split}
\end{align}
Then \eqref{eq-4.13} and \eqref{eq-4.14}
guarantee that
$$
H_{1}(\partial_{z}f)= H_{1}(\partial_{\overline{z}}f)=\infty.
$$
These complete the proof of the third statement in the example.

(4)   Since
 \eqref{eq-4.13} and \eqref{eq-4.14} yield that
$$
H_{\infty}(\partial_{\overline{z}}f)
 =H_{\infty}(\partial_{z}f)
= \sup_{z\in \mathbb{D}}|\partial_{z}f(z)|
\geq H_{1}(\partial_{z}f),
$$
 we see that
  $\| D_{f}\|$  is unbounded in $\mathbb{D}$.
This and \eqref{eq-4.14} ensure that $f$ is not a $(K,K')$-elliptic mapping for any $K\geq1$ and $K'\geq 0$, and hence, the proof of the fourth statement in the example is complete.
\epf

\subsection{Proof of Theorem \ref{thm-1.2}}
Theorem \ref{thm-1.2} follows from Examples \ref{ex-4.1}\eqref{ex-4.1-4},  \ref{ex-4.2}\eqref{ex-4.2-3} and  \ref{ex-4.3}\eqref{ex-4.3-3}.
\qed

\subsection{Proof of Theorem \ref{thm-1.4}}
Theorem \ref{thm-1.4}
 follows from
Examples \ref{ex-4.1}\eqref{ex-4.1-5}  and   \ref{ex-4.2}\eqref{ex-4.2-4}.
\qed

 \section{Acknowledgments}
 We sincerely thank professor David Kalaj for his valuable comments.
The first author (Jiaolong Chen) is partially supported by NNSF of China (No. 12071121),
NSF of Hunan Province (No. 2018JJ3327) and the construct program of the key discipline in Hunan Province.
The second author (Shaolin Chen) is partially supported by NNSF of China (No.  12071116), Hunan Provincial Natural Science Foundation of China (No. 2022JJ10001),
and the Key Projects of Hunan Provincial Department of Education (No. 21A0429).

\end{document}